 \documentclass[preprint,12pt,3p,fleqn]{elsarticle}
\pdfoutput=1

\usepackage{amssymb}
\usepackage{amsthm}
\usepackage{mathtools}
\usepackage{amsmath}
\usepackage{multirow}
\usepackage{booktabs}
\usepackage{graphicx}
\usepackage{epstopdf}
\usepackage{pdflscape}
\usepackage{amsthm}
\usepackage{amsmath}
\usepackage[english]{babel}
\usepackage{subfigure}
\usepackage{threeparttable}

\usepackage{graphicx}
\usepackage{subfigure}

\usepackage{tabularx}
\usepackage{graphicx}
\usepackage{threeparttable}
\usepackage[figuresright]{rotating}
\usepackage{float}
\allowdisplaybreaks

\theoremstyle{definition}

\theoremstyle{remark}

\usepackage{hyperref}

\usepackage{xcolor}

\journal{XXXX}

\begin{document}

\begin{frontmatter}

\title{A class of new high-order finite-volume TENO schemes for hyperbolic conservation laws with unstructured meshes}

\author[a,b]{Zhe Ji}
\ead{jizhe@nwpu.edu.cn}

\author[c]{Tian Liang}
\ead{tliangae@connect.ust.hk}

\author[c,d,e,f]{Lin Fu\corref{cor1}}
\ead{linfu@ust.hk}
\cortext[cor1]{Corresponding author.}

\address[a]{School of Software, Northwestern Polytechnical University, Taicang, Jiangsu, China}
\address[b]{Yangtze River Delta Research Institute of NPU, Taicang, Jiangsu, China}

\address[c]{Department of Mechanical and Aerospace Engineering, The Hong Kong University of Science and Technology, Clear Water Bay, Kowloon, Hong Kong}
\address[d]{Department of Mathematics, The Hong Kong University of Science and Technology, Clear Water Bay, Kowloon, Hong Kong}
\address[e]{Center for Ocean Research in Hong Kong and Macau (CORE), The Hong Kong University of Science and Technology, Clear Water Bay, Kowloon, Hong Kong}
\address[f]{Shenzhen Research Institute, The Hong Kong University of Science and Technology, Shenzhen, China}


%
\begin{abstract}

    The recently proposed high-order TENO scheme [Fu et al., Journal of Computational Physics, 305(2016): 333-359] has shown great potential in predicting complex fluids owing to the novel weighting strategy, which ensures the high-order accuracy, the low numerical dissipation, and the sharp shock-capturing capability. However, the applications are still restricted to simple geometries with Cartesian or curvilinear meshes. In this work, a new class of high-order shock-capturing TENO schemes for unstructured meshes are proposed. Similar to the standard TENO schemes and some variants of WENO schemes, the candidate stencils include one large stencil and several small third-order stencils. Following a strong scale-separation procedure, a tailored novel ENO-like stencil selection strategy is proposed such that the high-order accuracy is restored in smooth regions by selecting the candidate reconstruction on the large stencil while the ENO property is enforced near discontinuities by adopting the candidate reconstruction from smooth small stencils. The nonsmooth stencils containing genuine discontinuities are explicitly excluded from the final reconstruction, leading to excellent numerical stability. Different from the WENO concept, such unique sharp stencil selection retains the low numerical dissipation without sacrificing the shock-capturing capability. The newly proposed framework enables arbitrarily high-order TENO reconstructions on unstructured meshes. For conceptual verification, the TENO schemes with third- to sixth-order accuracy are constructed. Without parameter tuning case by case, the performance of the proposed TENO schemes is demonstrated by examining a set of benchmark cases with broadband flow length scales.

\end{abstract}

\end{frontmatter}

\section{Introduction}

    For computational fluid dynamics (CFD), one core research topic is to develop high-order stable numerical methods for solving the governing hyperbolic conservation laws, i.e. the Euler or Navier–Stokes equations. The solution complexity of this set of initial-boundary value PDEs resides at that discontinuities may appear even when the initial condition is sufficiently smooth \cite{fu2021shock}\cite{griffin2021velocity}\cite{bai2022assessment}. The basic requirements for modern numerical methods are that the high-order accuracy is achieved for smooth flows and the solution monotonicity is preserved in the vicinity of discontinuities \cite{shu2009high}\cite{Pirozzoli2011}. For more complex fluid simulations, the low-dissipation property is also essential such that the small-scale flow structures can be resolved with high fidelity \cite{johnsen2010assessment}\cite{fu2019LES}\cite{fu2019improved}. 

    For flows over complex geometries, the established high-order numerical methods include, e.g., the conventional finite-volume \cite{shu2003high}\cite{ollivier2002high}\cite{diskin2007accuracy} method, the discontinuous Galerkin (DG) \cite{cockburn1989tvb}\cite{cockburn1998runge} method, and the flux-reconstruction (FR) \cite{huynh2007flux}\cite{witherden2016high} method. The finite-volume method typically relies on the data from multi-layer neighboring cells for the targeted high-order reconstruction. On the other hand, DG and FR methods can achieve high-order accuracy by deploying high-order polynomials within the mesh element. This compactness allows for a dramatic reduction of programming complexity and inter-node data transfer for high-performance parallel computing with unstructured meshes \cite{zhou2001numerical}. Therefore, for low-Mach flows, DG and FR methods are more and more popular for high-fidelity simulations, e.g. for large-eddy simulations (LES) \cite{flad2017use}\cite{frere2017application}\cite{krank2019multiscale} and direct numerical simulations (DNS) \cite{collis2002discontinuous}\cite{renac2015aghora}\cite{gempesaw2011multi}. 

    However, when high-Mach flows with shockwaves are concerned, the deployment of the shock-capturing concept in these frameworks becomes inevitable and challenging. For the DG and FR frameworks, the most popular choice is to introduce the artificial viscosity method \cite{abbassi2014shock}\cite{haga2019robust}\cite{vandenhoeck2019implicit}\cite{witherden2016high} for its simplicity in terms of implementation. However, the performance of the artificial viscosity method is typically unsatisfactory and the built-in parameters need intensive case-by-case tuning. One alternative avenue is to resort to a hybrid strategy, which first detects the shockwaves based on a certain indicator and locally deploys the finite-volume shock-capturing concept, e.g. the Total Variation Diminishing (TVD) \cite{van1979towards}\cite{van1977towards}\cite{van1974towards} method, the essentially non-oscillatory (ENO) \cite{Harten1987} method, the weighted ENO (WENO) \cite{Liu1994} method, and the Total Variation Bounding (TVB) \cite{shu1987tvb} method. {\color{black}Other relevant work includes the class of subcell finite volume limiters introduced for DG schemes with a posteriori MOOD detector for Cartesian and unstructured meshes \cite{boscheri2017arbitrary}\cite{boscheri2014direct}\cite{boscheri2014lagrangian}.}

    In terms of the TVD-based schemes for unstructured meshes, Barth and Jespersen \cite{barth1989design} propose a novel oscillation-free second-order scheme involving a limiter acting on a linear reconstruction polynomial. Venkatakrishnan \cite{venkatakrishnan1995convergence} further improves it by making it differentiable and thus better convergence property for steady state solutions is obtained. To extend these limiters with high-order accuracy for smooth flows, Li and Ren \cite{li2012multi} propose the so-called weighted biased averaging procedure (WBAP) limiter by employing non-linear weighting coefficients for the reconstruction and up to fourth-order accuracy on unstructured meshes is achieved. Michalak and Ollivier-Gooch \cite{michalak2009accuracy} develop the accuracy preserving limiter for high-order accurate solution of the Euler equations, where fourth-order accuracy is achieved with minimal sacrifices in convergence properties comparing to existing second-order approaches. 

    However, when concerning even higher-order reconstruction, the further development of limiter-based shock-capturing schemes becomes incredibly cumbersome. Instead, the WENO paradigm \cite{Liu1994}, originating from the ENO concept \cite{Harten1987}, is able to provide arbitrarily high-order reconstruction and sharp shock-capturing capability within a unified framework. Different from ENO, WENO schemes utilize a weighted average of several low-order candidates to form the high-order reconstruction on the full stencil. The nonlinear weights are computed dynamically based on the smoothness indicators of candidate stencils such that the optimal accuracy order is restored asymptotically in smooth regions and the ENO property is enforced near discontinuities. While great success has been achieved for applications on structured meshes with WENO schemes, it is, however, non-trivial to extend the WENO weighting procedure to unstructured meshes. The main issue is that, for some mesh topologies, the optimal linear weights may become negative inducing numerical instabilities or even do not exist \cite{hu1999weighted}\cite{zhang2009third}\cite{shi2002technique}. One way to solve this issue is to decrease the targeted accuracy order on the combined large stencil to that of the small candidate stencils \cite{cheng2008third}\cite{dumbser2007arbitrary}\cite{dumbser2007quadrature}\cite{liu2013robust}. Alternatively, the central WENO (CWENO) schemes \cite{levy1999central}\cite{levy2000compact}\cite{capdeville2008central}\cite{cravero2018cweno}, which feature stencils of different sizes and tailored weighting strategy, can avoid these numerical difficulties since the choice of the linear weights can be arbitrary as long as their total sum equals one. The latest variants of CWENO schemes include the multi-resolution WENO schemes \cite{zhu2019new}\cite{zhu2020new} and the WENO schemes with adaptive orders \cite{balsara2020efficient}. {\color{black} The CWENO/CWENOZ schemes have been extended to mixed-element unstructured meshes and multicomponent flows using unstructured meshes \cite{tsoutsanis2021arbitrary}\cite{tsoutsanis2021cweno}.} {\color{black}The WENO adaptive-order scheme has also been recently used for general PNPM schemes (reconstructed DG) on fixed and moving meshes \cite{boscheri2019high}.} 

    More recently, Fu et al. \cite{fu2016family}\cite{fu2017targeted}\cite{fu2018new}\cite{fu2021very}\cite{takagi2021novel} propose a family of high-order targeted ENO (TENO) schemes for solving the hyperbolic conservation laws. By introducing the novel ENO-like stencil selection strategy, the TENO schemes feature low numerical dissipation and sharp shock-capturing property. Compared to WENO schemes, one unique property of TENO is that the background linear scheme can be restored exactly for up to intermediate wavenumbers without sacrificing the robustness for capturing strong discontinuities. In this way, the performance of TENO schemes can be well controlled by optimizing the spectral properties of the underlying background linear schemes. The TENO-family schemes have been widely deployed to different complex fluid simulations on structured meshes, e.g. the multi-phase flows \cite{haimovich2017numerical}, the detonation simulations \cite{dongdetonation}, the magnetohydrodynamics (MHD) \cite{fu2019high}\cite{fu2022efficient}, the compressible gas dynamics \cite{sun2016boundary}\cite{zhao2020104439}\cite{zhang2020simple}\cite{zhang2019towards}\cite{fardipour2020development}\cite{fu2019hybrid}, the turbulent flows \cite{di2020htr}\cite{motheau2020capturing}\cite{Lusher2020}\cite{lefieux2019dns}\cite{lusher2019assessment}, etc. 

    As far as the authors' knowledge, the high-order TENO schemes have not been extended to the unstructured meshes, and consequently the applications are restricted to canonical simulations. In this work, a new class of high-order TENO schemes on unstructured meshes are proposed. Similar to the idea in the standard TENO schemes \cite{fu2016family}\cite{fu2017targeted}\cite{fu2018new} and the WENO variants \cite{levy1999central}\cite{capdeville2008central} \cite{zhu2016new}\cite{balsara2016efficient}, the candidate stencils involve one large central-biased stencil and several small directional stencils. A tailored novel ENO-like stencil selection strategy is proposed, which ensures that the final scheme is equipped with the high-order reconstruction on the large stencil in smooth regions and with the low-order reconstruction on small smooth stencils near the discontinuities. Except for the constant $C_T$ in the weighting strategy, which separates the discontinuities from smooth scales in spectral space and can be determined by aprior analyses, no additional case-by-case empirical parameters are involved in the proposed TENO schemes. Moreover, arbitrarily high-order accurate TENO schemes can be constructed with this unified framework. For the conceptual demonstration, the TENO schemes with third- to sixth-order accuracy are constructed and validated with a set of critical benchmark simulations.

    The remainder of this paper is organized as follows. (1) In section 2, the high-order finite-volume framework for solving the hyperbolic conservation laws is briefly reviewed as well as the WENO and CWENO reconstruction schemes; (2) In section 3, the new high-order TENO schemes for unstructured meshes are developed; (3) In section 4, the numerical method for flux evaluations is discussed;
    (4) In section 5, the performance of the proposed low-dissipation TENO schemes is demonstrated by conducting a wide range of critical benchmark simulations;
    (5) Conclusions and remarks are given in the last section.

\section{Fundamentals of the unstructured finite-volume methods}

    In this section, the basic concepts of the unstructured finite-volume method including the classical high-order WENO and CWENO reconstruction schemes will be outlined.
    \subsection{Concepts of finite-volume method}
        For three-dimensional unsteady Euler equations, as typical hyperbolic conservation laws, the conservative form can be written as
        \begin{equation}
            \label{eq:fv_governing_equation}
            \frac{\partial\textbf{U}}{\partial t}+\nabla\cdot(\textbf{F}(\textbf{U})) = 0,
        \end{equation}
        where $\textbf{U}=\textbf{U}(\textbf{x},t)$ denotes the conservative variables and $\textbf{F}(\textbf{U})=(\textbf{F}^x(\textbf{U}), \textbf{F}^y(\textbf{U}), \textbf{F}^z(\textbf{U}))$ the nonlinear convection flux functions. Assuming that the computational domain $\bf{\Omega}$ is partitioned by elements of various shapes, e.g. triangle, quadrilateral in 2D, or tetrahedral, hexahedral in 3D, integrating Eq.~(\ref{eq:fv_governing_equation}) over one control element $\mathbb{I}_i$ results in the following semi-discrete form
        \begin{equation}
            \label{eq:fv_governing_equation_semi_discrete}
            \frac{d{\overline {\textbf{U}}}_{\mathbb{I}_i}} {dt} = -\dfrac{1}{|V_{\mathbb{I}_i}|}\sum^{N_f}_{j}\sum^{N_q}_{\kappa}F^{\textbf{n}_{\mathbb{I}_{i,j}}}\left(\textbf{U}^n_{\mathbb{I}_{i,j,L}}(\textbf{x}_{\mathbb{I}_{i,j,\kappa}},t),\textbf{U}^n_{\mathbb{I}_{i,j,R}}(\textbf{x}_{\mathbb{I}_{i,j,\kappa}},t)\right)\omega_{\kappa}|\textbf{A}_{\mathbb{I}_{i,j}}|,
        \end{equation}
        where $\overline {\textbf{U}}_{\mathbb{I}_i}$ is the volume-averaged conservative variable, $V_{\mathbb{I}_i}$ the volume of control element ${\mathbb{I}_i}$, $N_f$ the cell face number of ${\mathbb{I}_i}$, $N_q$ the number of quadrature points deployed for the high-order surface integral approximation, $F^{\textbf{n}_{\mathbb{I}_{i,j}}}$ the numerical flux in the normal direction (pointing outwards) of face $j$, $\textbf{U}^n_{\mathbb{I}_{i,j,L}}$ and $\textbf{U}^n_{\mathbb{I}_{i,j,R}}$ the left- and right-biased approximated solutions on the interface respectively, $\omega_{\kappa}$ the weight for Gaussian integration point $\textbf{x}_{\mathbb{I}_{i,j,\kappa}}$, and $|\textbf{A}_{\mathbb{I}_{i,j}}|$ the surface area of face $j$.
        
        In order to advance the cell averages of the solution in time, the third-order strong-stability-preserving (SSP) Runge-Kutta method \cite{gottlieb2001strong} is deployed for all the following investigations unless explicitly specified. The explicit formulas are given as
        \begin{equation}
                 \label{eq:fv_time_discretization_01}
                 \begin{array}{l}
                 {\overline {\textbf{U}}}^{(1)}_{\mathbb{I}_i}={\overline {\textbf{U}}}^{n}_{\mathbb{I}_i}+\Delta t\mathcal{R}({\overline {\textbf{U}}}^{n}_{\mathbb{I}_i}),\\
                 {\overline {\textbf{U}}}^{(2)}_{\mathbb{I}_i}=\dfrac{3}{4}{\overline {\textbf{U}}}^{n}_{\mathbb{I}_i}+\dfrac{1}{4}{\overline {\textbf{U}}}^{(1)}_{\mathbb{I}_i}+\dfrac{1}{4}\Delta t\mathcal{R}({\overline {\textbf{U}}}^{(1)}_{\mathbb{I}_i}),\\
                 {\overline {\textbf{U}}}^{n+1}_{\mathbb{I}_i}=\dfrac{1}{3}{\overline {\textbf{U}}}^{n}_{\mathbb{I}_i}+\dfrac{2}{3}{\overline {\textbf{U}}}^{(2)}_{\mathbb{I}_i}+\dfrac{2}{3}\Delta t\mathcal{R}({\overline {\textbf{U}}}^{(2)}_{\mathbb{I}_i}),
                 \end{array}
         \end{equation}
        where $\Delta t$ is the timestep size, and $\mathcal{R}({\overline {\textbf{U}}}^{n}_{\mathbb{I}_i})$ denotes the right hand side of Eq.~(\ref{eq:fv_governing_equation_semi_discrete}).
        
        Based on the standard finite-volume framework, the remaining numerical issue is to develop the high-order scheme to reconstruct the cell interface data, i.e. $\textbf{U}^n_{\mathbb{I}_{i,j,L}}(\textbf{x}_{\mathbb{I}_{i,j,\kappa}},t)$ and $\textbf{U}^n_{\mathbb{I}_{i,j,R}}(\textbf{x}_{\mathbb{I}_{i,j,\kappa}},t)$, based on the known cell-averaged solutions, and to evaluate the cell interface flux $F^{\textbf{n}_{\mathbb{I}_{i,j}}}\left(\textbf{U}^n_{\mathbb{I}_{i,j,L}}(\textbf{x}_{\mathbb{I}_{i,j,\kappa}},t),\textbf{U}^n_{\mathbb{I}_{i,j,R}}(\textbf{x}_{\mathbb{I}_{i,j,\kappa}},t)\right)$.

    \subsection{The high-order linear reconstruction schemes}

        Given a polynomial $\mathcal{P}_{\mathbb{I}_i}(x,y,z)$ of order $r$ in cell ${\mathbb{I}_i}$, if the reconstructed cell-averaged value equals ${\overline {\textbf{U}}}_{\mathbb{I}_i}$, i.e.
        \begin{equation}
            \label{eq:reconstruction_01}
            {\overline {\textbf{U}}}_{\mathbb{I}_i} = \dfrac{1}{|V_{\mathbb{I}_i}|}\int_{V_{\mathbb{I}_i}}\textbf{U}(x,y,z)dV = \dfrac{1}{|V_{\mathbb{I}_i}|}\int_{V_{\mathbb{I}_i}}\mathcal{P}_{\mathbb{I}_i}(x,y,z)dV,
        \end{equation}
        then $r+1$ order of accuracy is achieved. For unstructured meshes involving elements of various shapes, it is beneficial to transform all the elements from physical space ($\textbf{X}=(x,y,z)$) to a reference space ($\mathbf{\Xi}=(\xi,\eta,\zeta)$) in order to reduce the scaling effects \cite{dumbser2007quadrature}\cite{dumbser2007arbitrary}. Following \cite{tsoutsanis2011weno}\cite{tsoutsanis2014weno}, all types of elements are first decomposed into basic elements of triangular or tetrahedral shape depending on the dimension of the problem. Then the coordinates of one of the new basic element are used as the reference to perform the transformation. Note that the spatial transformation does not change the cell-averaged solutions, i.e.
        \begin{equation}
            \label{eq:reconstruction_02}
            {\overline {\textbf{U}}}_{\mathbb{I}_i} = \dfrac{1}{|V_{\mathbb{I}_i}|}\int_{V_{\mathbb{I}_i}}\textbf{U}(x,y,z)dV \equiv \dfrac{1}{|V^{'}_{\mathbb{I}_i}|}\int_{V^{'}_{\mathbb{I}_i}}\textbf{U}(\xi,\eta,\zeta)d\xi d\eta d\zeta.
        \end{equation}
        For detailed descriptions, one is referred to \cite{tsoutsanis2011weno}\cite{tsoutsanis2019stencil}.
        
        To calculate reconstructed values in the targeted cell ${\mathbb{I}_0}$, a compact stencil $\mathcal{S}$ is built with neighboring cells of ${\mathbb{I}_0}$. The amount of elements in $\mathcal{S}$, i.e. $N_s$, is determined by the total number of polynomial coefficients, $N_k$, which is calculated by
        \begin{equation}
            \label{eq:reconstruction_03}
            N_k=\dfrac{1}{n_d!}\prod_{m=1}^{n_d}(r+m),
        \end{equation}
        where $n_d$ is the number of dimension. In order to maintain good numerical stability with an acceptable computational cost, it is recommended to choose $N_s=2\cdot N_k$ according to various studies \cite{tsoutsanis2014weno}\cite{tsoutsanis2019stencil}\cite{tsoutsanis2018improvement}, and we follow this setup for all the investigations in the present work.
        
        The reconstruction polynomial $\mathcal{P}(\xi,\eta,\zeta)$ can be expanded in $\mathcal{S}$ over the local polynomial basis functions $\psi_l(\xi,\eta,\zeta)$ as
        \begin{equation}
            \label{eq:reconstruction_04}
            \mathcal{P}(\xi,\eta,\zeta)=\sum_{l=0}^{N_k}a_l\psi_l(\xi,\eta,\zeta)=\overline {\textbf{U}}_0+\sum_{l=1}^{N_k}a_l\psi_l(\xi,\eta,\zeta),
        \end{equation}
        where $\overline {\textbf{U}}_0$ denotes the cell-averaged solution of the target cell ${\mathbb{I}_0}$, $a_l$ the degrees of freedom to be calculated later. The basis functions $\psi_l(\xi,\eta,\zeta)$ are defined as
        \begin{equation}
            \label{eq:reconstruction_06}
            \begin{array}{c}
            \psi_l(\xi,\eta,\zeta) \equiv \varphi_l(\xi,\eta,\zeta)-\dfrac{1}{|V^{'}_{0}|}\int_{V^{'}_{0}}\varphi_ld\xi d\eta d\zeta \text{, } l=1,2,..,N_k,\\
            \end{array}
        \end{equation}
        where
        \begin{equation}
            \label{eq:reconstruction_06_222}
            \begin{array}{c}
            {\varphi_l}=\xi,\text{ }\eta,\text{ }\zeta,\text{ }\xi^2,\text{ }\eta^2,\text{ }\zeta^2,\text{ } \xi\cdot\eta,...,
            \end{array}
        \end{equation}
        to ensure that the constraint of Eq.~(\ref{eq:reconstruction_01}) for the targeted cell ${\mathbb{I}_0}$ is satisfied irrespective of choices of $a_l$.
        
        Based on the constraint that the reconstructed cell-averaged values for all of the cells ${\mathbb{I}_s}$ in $\mathcal{S}$ should be equal to the corresponding cell-averaged solutions, one can obtain
        \begin{equation}
            \label{eq:reconstruction_05}
            \begin{array}{c}
            \int_{V^{'}_{s}}\mathcal{P}(\xi,\eta,\zeta)d\xi d\eta d\zeta=|V^{'}_{s}|\overline {\textbf{U}}_0+\sum_{l=1}^{N_k}\int_{V^{'}_{s}}a_l\psi_l(\xi,\eta,\zeta)d\xi d\eta d\zeta = |V^{'}_{s}|\overline {\textbf{U}}_{s},\\
            s = 1,...,N_s.
            \end{array}
        \end{equation}
        To find the unknown degrees of freedom $a_l$, Eq.~(\ref{eq:reconstruction_05}) can be rewritten into
        \begin{equation}
            \label{eq:reconstruction_07}
            \sum_{l=1}^{N_k}\mathcal{A}_{sl}a_l=b_s,\text{ }s = 1,...,N_s,
        \end{equation}
        where
        \begin{equation}
            \label{eq:reconstruction_08}
            \mathcal{A}_{sl}=\int_{V^{'}_{s}}\psi_l(\xi,\eta,\zeta)d\xi d\eta d\zeta,\text{ and }b_s=|V^{'}_{s}|(\overline{\textbf{U}}_s-\overline{\textbf{U}}_0).
        \end{equation}
        The matrix $\mathcal{A}_{sl}$ contains only the geometry information of each element in the considered stencil $\mathcal{S}$, and hence can be pre-calculated before the simulation. Since the matrix $\mathcal{A}^{T}_{ls}\mathcal{A}_{sl}$ is invertible, $a_l$ is obtained as
        \begin{equation}
            \label{eq:reconstruction_09}
            a_l=(\mathcal{A}^{T}_{ls}\mathcal{A}_{sl})^{-1}\mathcal{A}^{T}_{ls}b_s=\widehat{\mathcal{A}}_{ls}b_s,
        \end{equation}
        where $\widehat{\mathcal{A}}_{ls}$ can be further simplified by deploying a QR decomposition based on the Householder transformation \cite{tsoutsanis2019stencil}\cite{tsoutsanis2018improvement} as
        \begin{equation}
            \label{eq:reconstruction_10}
            \widehat{\mathcal{A}}_{ls}=\left((QR)^T(QR)\right)^{-1}\mathcal{A}^{T}_{ls}=\left(R^TR\right)^{-1}\mathcal{A}^{T}_{ls}.
        \end{equation}
    \subsection{The high-order WENO reconstruction schemes}

    For a hyperbolic system, the above high-order linear reconstruction schemes are not suitable for handling discontinuities, such as shockwaves, presented in the computational domain. As in classical shock-capturing WENO schemes \cite{Jiang1996}, the main idea is to obtain the high-order oscillation-free reconstruction based on a nonlinear convex combination of several candidate stencils. The discontinuity can be captured sharply and stably by enforcing the ENO property through a nonlinear adaptation between the candidate reconstructions.
    
    However, depending on the mesh topologies, the optimal linear weights may not exist for some circumstances \cite{hu1999weighted}\cite{zhang2009third}\cite{shi2002technique}. To ensure the numerical robustness, the WENO scheme in \cite{cheng2008third}\cite{dumbser2007arbitrary}\cite{dumbser2007quadrature}\cite{liu2013robust} is considered, with which the targeted accuracy order of the combined reconstruction is identical to those of the candidate stencils.
    
    Specifically, based on a set of candidate stencils ($\mathcal{S}_0$,$\mathcal{S}_1$, \dots, $\mathcal{S}_K$) with linear reconstructions $\mathcal{P}_k(\xi, \eta, \zeta), k=0,\dots,K$ featuring the same targeted accuracy order, the final nonlinearly combined WENO reconstruction is given as

        \begin{equation}
        \label{eq:WENO_reconstruction}    
        \mathcal{P}^{\rm{W}}(\xi, \eta, \zeta) =  \sum\limits_{k = 0}^{K} {w_k} \mathcal{P}_k(\xi, \eta, \zeta) ,
        \end{equation}
    where the nonlinear weight $w_k$ is defined as 
        \begin{equation}
        \label{eq:WENOweight}
        \omega_{k} = \frac{\alpha_k}{\sum_{k=0}^{K}\alpha_k},  \text{  and    } \alpha_k = \frac{d_k}{{({\beta _k} + \varepsilon )}^4},  \text{ } k = 0,\cdots,K,
        \end{equation}
    with $\varepsilon = 10^{-14}$ and the smoothness indicators measured following the original definition in \cite{Jiang1996} as
        \begin{equation}
        \label{eq:WENO_smoothness}
        {\beta_k} =  \sum\limits_{q = 1}^{r} \oint_{V^{'}_{\mathbb{I}_0}} (\mathcal{D}^q \mathcal{P}_k(\xi, \eta, \zeta))^2 \,(d\xi, d\eta, d\zeta).
        \end{equation}
    Since the optimal linear weights $d_k$ are not designed to generate an even higher-order reconstruction scheme on the combined stencil, the values can be assigned flexibly. Previous investigations \cite{cheng2008third}\cite{dumbser2007arbitrary}\cite{dumbser2007quadrature} reveal that the performance of the resulting WENO scheme can be improved by increasing the weight of the central-biased candidate stencil, which is less dissipative by construction. In practice, we set $d_K = 10^4$ for the central-biased stencil $\mathcal{S}_K$ and $d_k = 1, k=0,\dots,K-1$ for the other directional stencils ($\mathcal{S}_0$,$\mathcal{S}_1$, \dots, $\mathcal{S}_{K-1}$), empirically.
    
    Note that, all the linear candidate reconstructions $\mathcal{P}_k(\xi, \eta, \zeta)$ of particular accuracy orders are constructed in the way described in the previous subsection, and satisfy the constraints of matching the cell-averaged solutions of involved stencil cells. All of them are obtained by solving the overdetermined linear systems with the same constrained least-squares technique.
    Although there is a debate on which variable to be reconstructed by the WENO-type schemes, it is generally believed that the numerical oscillations can be better suppressed when the WENO reconstruction is carried out for characteristic variables than primitive or conservative variables \cite{Jiang1996}\cite{fu2016family}. 

    \subsection{The high-order CWENO reconstruction schemes}

        In CWENO schemes \cite{levy1999central}\cite{levy2000compact}\cite{capdeville2008central}\cite{cravero2018cweno}, the candidate stencils include one large central-biased stencil $\mathcal{S}_K$ and several small directional stencils ($\mathcal{S}_0$,$\mathcal{S}_1$, \dots, $\mathcal{S}_{K-1}$), with which the high-order reconstruction scheme $\mathcal{P}_K(\xi, \eta, \zeta)$ and the low-order reconstruction schemes $\mathcal{P}_k(\xi, \eta, \zeta), k=0,\dots,K-1$ are constructed, respectively. Note that, for good numerical robustness, the third-order reconstructions are adopted for small stencils. The most important ingredient of CWENO schemes is the introduction of an optimal reconstruction scheme as
        \begin{equation}
        \label{eq:CWENOoptimal}  
        \mathcal{P}^{\rm{opt}}(\xi, \eta, \zeta) = \frac{1}{d_K} \left(\mathcal{P}_K(\xi, \eta, \zeta) - \sum\limits_{k = 0}^{K-1} {d_k} \mathcal{P}_k(\xi, \eta, \zeta)\right),
        \end{equation}
        where the contributions of the lower-order reconstruction schemes are subtracted from the high-order reconstruction. Subsequently, the CWENO nonlinear scheme is given as a convex combination of the optimal candidate reconstruction scheme and other low-order candidate reconstruction schemes following
        \begin{equation}
        \label{eq:CWENO_reconstruction}    
        \mathcal{P}^{\rm{CW}}(\xi, \eta, \zeta) = \frac{w_K}{d_K} \left(\mathcal{P}_K(\xi, \eta, \zeta) - \sum\limits_{k = 0}^{K-1} {d_k} \mathcal{P}_k(\xi, \eta, \zeta)\right) + \sum\limits_{k = 0}^{K-1} {w_k} \mathcal{P}_k(\xi, \eta, \zeta),
        \end{equation}
        where $\omega_{k}$ denotes the nonlinear weights and can be computed in a similar way as WENO schemes by
        \begin{equation}
        \label{eq:CWENOweight}
        \omega_{k} = \frac{\alpha_k}{\sum_{k=0}^{K}\alpha_k},  \text{  and    } \alpha_k = \frac{d_k}{{({\beta _k} + \varepsilon )}^4}, \text{ } k = 0,\cdots,K,
        \end{equation}
       where $\varepsilon = 10^{-14}$. The smoothness indicator $\beta _k$ is defined to be the same as that in WENO schemes. Unlike the WENO schemes, the linear weights $d_k$ should satisfy
        \begin{equation}
        \label{eq:CWENOOptimalLinearweight}
        \sum_{k=0}^{K}d_k=1,\text{  and   }d_K=1-\dfrac{1}{d_K^{'}},
        \end{equation}
       where $d_K^{'}$ can be taken arbitrarily and $d_k,\text{ }k=0,\dots,K-1$ is calculated by
        \begin{equation}
        \label{eq:CWENOOptimalLinearweight02}
        d_k=\dfrac{1-d_K}{K}.
        \end{equation}
        When the local flow scale is smooth, $\omega_{k} \approx d_k$ and the high-order reconstruction is restored asymptotically with $\mathcal{P}^{\rm{CW}}(\xi, \eta, \zeta) \approx \mathcal{P}_K(\xi, \eta, \zeta)$; when approaching discontinuities, the contribution from the large candidate stencil will vanish and the smooth small stencils will dominate the final reconstruction. Therefore, the overall accuracy order of the resulting CWENO scheme is determined by the reconstruction scheme on the large stencil and the corresponding performance in nonsmooth regions is governed by the reconstruction schemes on the small stencils. 
        
        Increasing the optimal linear weight for large stencil can improve the wave-resolution property of the resulting CWENO scheme since the final scheme is biased to the high-order reconstruction. However, a too large value of $d_K$ may also degrade the shock-capturing capability. In practice, we set $d_K^{'} = 10^4$, i.e., $d_K=1-{10^{-4}}$, for the large stencil and $d_k = \frac{10^{-4}}{K}, k=0,\dots,K-1$ for the small directional stencils, empirically. {\color{black}Note that different $d_K^{'}$ values can be adopted for CWENO schemes, i.e., $10^3$, $10^7$, and $10^{15}$ \cite{tsoutsanis2021arbitrary}. And one single $d_K^{'}$ value may not be suitable for all the discretization orders. Considering the fact that all the present TENO schemes of different accuracy orders employ the same set of parameters without tuning case by case, for fair comparisons, we adopt this uniform optimal $d_K^{'}$ value for CWENO schemes as well in the present study if not mentioned otherwise.}
        
        It is also worth noting that, to achieve the same accuracy order, the cost and complexity of the CWENO scheme are much smaller than the WENO scheme, since all the candidate stencils of WENO feature the same reconstruction order as the targeted order of the resulting nonlinear scheme. Also, the compactness of the small directional candidate stencils in CWENO increases the chance that at least one candidate is smooth in the high-gradient nonsmooth region and therefore enhances the numerical robustness when compared to WENO schemes. {\color{black} The CWENO/CWENOZ schemes have been extended to mixed-element unstructured meshes and multicomponent flows using unstructured meshes \cite{tsoutsanis2021arbitrary}\cite{tsoutsanis2021cweno}.}
\section{The new class of high-order unstructured TENO reconstruction schemes}

    In this section, a new class of high-order unstructured TENO reconstruction schemes will be developed. The unstructured TENO weighting strategy including the candidate stencil arrangement, the strong scale separation, the novel ENO-type stencil selection will be elaborated in detail. 

    \subsection{The reconstruction candidate stencils}

        Similar to previous TENO \cite{fu2019low} and CWENO \cite{levy1999central} schemes, the core of the newly proposed unstructured TENO scheme is to design a set of candidate stencils involving both the large central-biased stencil and several small directional stencils. The desired overall high-order reconstruction scheme will be constructed based on the large central-biased stencil $\mathcal{S}_K$ for resolving smooth scales, whereas a set of third-order reconstruction schemes will be constructed with the small directional stencils ($\mathcal{S}_0$,$\mathcal{S}_1$, \dots, $\mathcal{S}_{K-1}$) for capturing discontinuities.

\begin{figure}[H]
            \centering
            \includegraphics[width=1.0\textwidth,trim=2 2 2 2,clip]{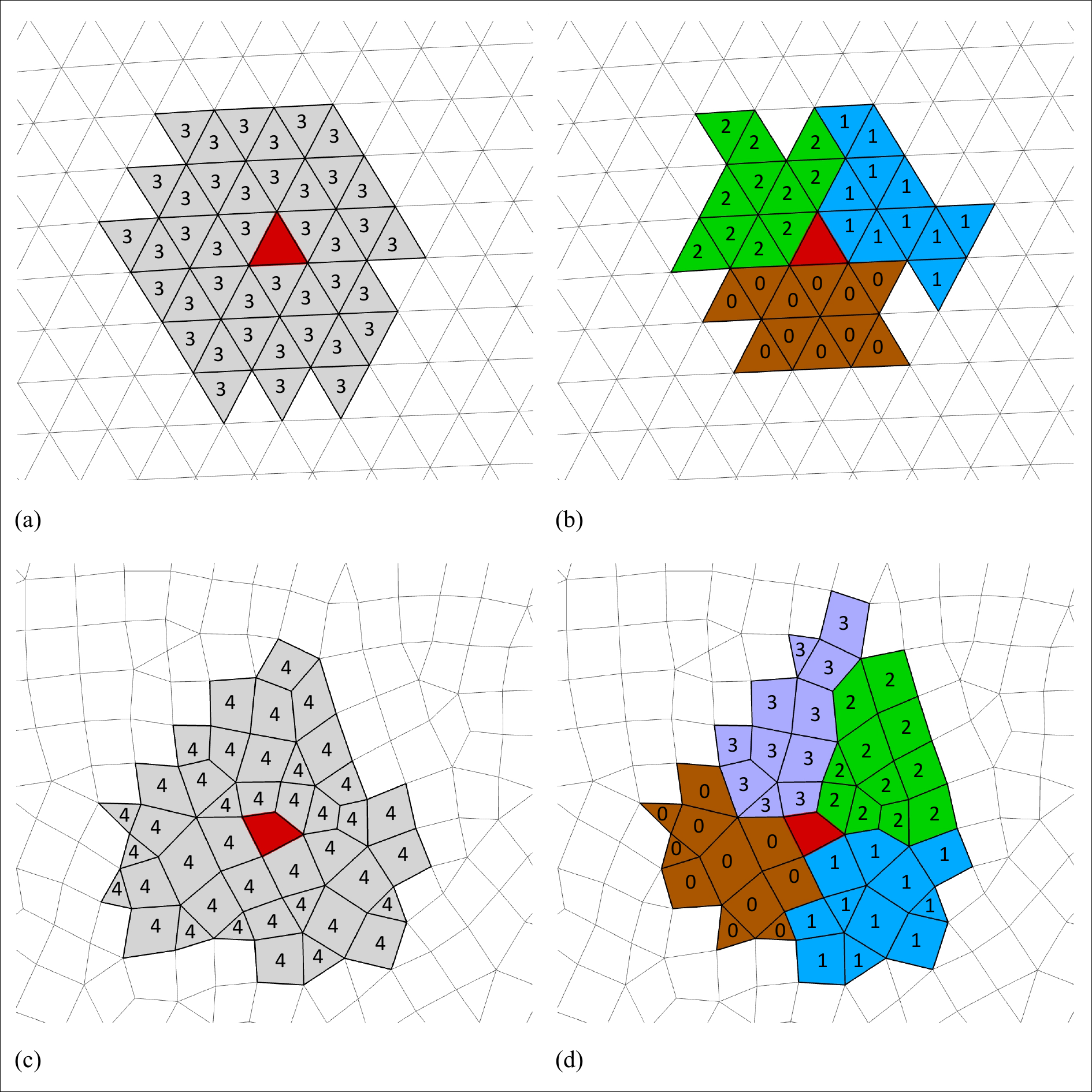}
            \caption{Sketch of (a) the sixth-order large central-biased candidate stencil $\mathcal{S}_3$ and (b) the third-order small directional candidate stencils ($\mathcal{S}_0$, $\mathcal{S}_1$, $\mathcal{S}_2$) for the TENO6 scheme on a triangular mesh; sketch of (c) the sixth-order large central-biased candidate stencil $\mathcal{S}_4$ and (d) the third-order small directional candidate stencils ($\mathcal{S}_0$, $\mathcal{S}_1$, $\mathcal{S}_2$, $\mathcal{S}_3$) for the TENO6 scheme on a mixed mesh. The index of each stencil is labelled in the figure and different stencils are filled with different colors. The cell considered for reconstruction $\mathbb{I}_0$ is marked with red color.}
            \label{Fig:stencil}
        \end{figure}

        Fig. \ref{Fig:stencil} shows the sketch of the central-biased and directional candidate stencils for the TENO6 scheme on different mesh topologies. In the present work, the central-biased stencil is built based on the Naive Cell based Algorithm (NCB) while the Type3 algorithm is employed to construct directional stencils following \cite{tsoutsanis2019stencil}. The stencil arrangement for unstructured meshes with other-type mesh elements follows a similar principle.

    \subsection{Strong scale separation}

        Similar to that in the Cartesian TENO schemes \cite{fu2019very}, in order to isolate discontinuities from the smooth scales, one core algorithm of the TENO weighting strategy is to measure the smoothness of candidate stencils with strong scale separation as 
        \begin{equation}
        \label{eq:new_scale_separation}
        {\gamma _k}  = \frac{1}{{({\beta _k} + \varepsilon )}^6},  \text{ } k = 0,\cdots,K,
        \end{equation}
        where $\varepsilon = 10^{-12}$ is a small value to avoid the zero denominator and $\beta _k$ denotes the smoothness function of candidate stencil $k$ following Eq.~(\ref{eq:WENO_smoothness}) \cite{Jiang1996}. {\color{black}Compared to the classical WENO \cite{Jiang1996} and CWENO \cite{levy1999central}\cite{levy2000compact}\cite{capdeville2008central} schemes, the exponent parameter is much larger for better scale separation and the low-dissipation property (which is difficult to achieve for WENO with this formula) will be retained by the following ENO-type stencil selection strategy. Note that, while the WENO-Z \cite{Borges2008}\cite{tsoutsanis2021cweno}\cite{tsoutsanis2021arbitrary} type measurement can be applied here, it is not adopted for reducing the computational costs in the context of unstructured meshes and the previous investigation with the Cartesian TENO schemes \cite{fu2019very} reveals that a 
       good scale-separation capability can be achieved with the present formula.}

    \subsection{The new ENO-like stencil selection}

        In contrast to classical WENO schemes \cite{Jiang1996}, where the contributions of the candidate stencils are combined in a smooth manner, for TENO, each stencil is explicitly judged to be a smooth or nonsmooth candidate by the so-called ENO-type stencil selection strategy. First, the measured smoothness indicator ${\chi}_k$ is renormalized as
        \begin{equation}
        \label{eq:normalize}
        {{\chi}_k } = \frac{\gamma _k}{{\sum\nolimits_{k = 0}^{K} {{\gamma _k}} }},
        \end{equation}
        and then subjected to a sharp cut-off function
        \begin{equation}
        \label{eq:cutoff}
        {\delta _k} = \left\{ {\begin{array}{*{20}{c}}
        0, &{\text{if }{\chi}_k  < {C_T},}\\
        1, &{\text{otherwise,}}
        \end{array}} \right.
        \end{equation}
        where the parameter $C_T$ determines the wavenumber, which separates the smooth and nonsmooth scales. As shown in \cite{fu2016family}\cite{fu2019LES}, the parameter $C_T$ can be either a constant or adaptively varying. If chosen as a constant, it is typically determined by the spectral analysis and by retaining the numerical robustness for high-Mach flows. In practice, the choice of $C_T$ is rather flexible between $10^{-5}$ and $10^{-7}$ \cite{fu2016family}, and the numerical robustness and the low-dissipation property can be demonstrated by conducting critical simulations with broadband flow scales without tuning parameter case-by-case \cite{fu2016family}\cite{fu2017targeted}\cite{fu2018new}, as will be shown in the numerical validation part of this work.
        
        With the present stencil arrangement, the large central-biased candidate stencil $\mathcal{S}_K$ features the highest accuracy order and the best spectral property and therefore is suitable for resolving the smooth flow scales. On the other hand, the set of small directional stencils is designed for capturing the possible discontinuity by a nonlinear adaptation among them. Specifically, the final data reconstruction on the cell interface can be determined in two steps as follows.
        
       (i) Group all the candidate stencils as $(\mathcal{S}_0, \mathcal{S}_1, \mathcal{S}_2, \dots, \mathcal{S}_K)$ and perform the ENO-type stencil selection following Eq.~(\ref{eq:normalize}) and Eq.~(\ref{eq:cutoff}). If the largest candidate stencil $\mathcal{S}_K$ is judged to be smooth, then the final reconstruction scheme is given as
        \begin{equation}
        \label{eq:TENO_reconstruction111}       
           \mathcal{P}^{\rm{T}}(\xi, \eta, \zeta) = \mathcal{P}_K (\xi, \eta, \zeta) ;
        \end{equation}

       (ii) If the largest candidate stencil $\mathcal{S}_k$ is judged to be nonsmooth, then the nonlinear adaptation among the small stencils will be activated to capture the discontinuities by enforcing the ENO property, i.e. group all the small candidate stencils as $(\mathcal{S}_0, \mathcal{S}_1, \mathcal{S}_2, \dots, \mathcal{S}_{K-1})$ and perform the ENO-type stencil selection. The final reconstruction scheme is given as 
        \begin{equation}
        \label{eq:TENO_reconstruction222}
        \mathcal{P}^{\rm{T}}(\xi, \eta, \zeta) =  \sum\limits_{k = 0}^{K-1} {w_k} \mathcal{P}_k(\xi, \eta, \zeta), 
        \end{equation}
        where the nonlinear weight 
        \begin{equation}
        \label{eq:TENO_weight}        
        {w_k} = \frac{{\delta _k}}{{\sum\nolimits_{k = 0}^{K-1}{\delta _k} }}, k = 0, 1, \dots, K-1 .
        \end{equation}
        Note that, for both stencil selection steps, the parameter $C_T$ is chosen to be the same, i.e., $C_T = 10^{-6}$. Further optimization of the choice of $C_T$ is beyond the scope of this paper and interested readers are referred to the discussions in \cite{fu2016family}\cite{fu2017targeted}\cite{fu2019LES}.

    \subsection{Analysis and discussions}

        As demonstrated by previous investigations \cite{fu2016family}\cite{fu2017targeted}\cite{fu2018new}, the TENO weighting strategy can accurately identify the nonsmooth discontinuities from smooth flow scales. On the other hand, in smooth regions, the desirable high-order reconstruction will be restored exactly without any compromise. In other words, in low-wavenumber regime, the nonlinear TENO scheme degenerates to the optimal linear scheme on the large stencil $\mathcal{S}_K$. In terms of resolving the small-scale features, such a property renders TENO schemes superior over WENO, for which the nonlinear adaptation (i.e., the nonlinear dissipation) is added even for low-wavenumber smooth regions.
        
        Another notable property of the proposed unstructured TENO schemes is that neither the optimal linear weights $d_k$ in WENO nor the artificial weights $d_K^{'}$ in CWENO schemes are necessary to specify, suggesting the highly general applicability.

\section{Numerical flux evaluations}
    In this section, we discuss the numerical method for flux evaluations to complete the spatial discretization under the finite-volume framework.

    \subsection{Riemann problem definition}
    As proved in \cite{toro2013riemann}, based on the rotational invariance property, the following relation,
    \begin{equation}
    F^{\textbf{n}_{\mathbb{I}_{i,j}}}(\mathbf{U})=\mathbf{T}^{-1}_{\mathbb{I}_{i,j}} \mathbf{F}^x\left(\mathbf{T}_{\mathbb{I}_{i,j}} \mathbf{U}\right),
    \end{equation}
    where $\mathbf{T}_{\mathbb{I}_{i,j}}$ denotes the unique rotation matrix for local cell face $j$ of cell $i$ (note that $\mathbf{T}_{\mathbb{I}_{i,j}}$ is fixed for specific mesh topology), applies to the Euler equations, and the projected numerical flux can be further evaluated by
    \begin{equation}
    \mathbf{K}_{i j}=\sum^{N_q}_{\kappa}F^{\textbf{n}_{\mathbb{I}_{i,j}}}\left(\textbf{U}^n_{\mathbb{I}_{i,j,L}}(\textbf{x}_{\mathbb{I}_{i,j,\kappa}},t),\textbf{U}^n_{\mathbb{I}_{i,j,R}}(\textbf{x}_{\mathbb{I}_{i,j,\kappa}},t)\right)\omega_{\kappa}|\textbf{A}_{\mathbb{I}_{i,j}}|=\sum^{N_q}_{\kappa} \mathbf{T}^{-1}_{\mathbb{I}_{i,j}} \mathbf{F}^x\left(\widehat{\mathbf{U}}_{L}, \widehat{\mathbf{U}}_{R}\right) \omega_{\kappa}|\textbf{A}_{\mathbb{I}_{i,j}}|,
    \end{equation}
    where the rotated variables are defined as 
    \begin{equation}
    \widehat{\mathbf{U}}_{L}=\mathbf{T}_{\mathbb{I}_{i,j}} \textbf{U}^n_{\mathbb{I}_{i,j,L}}(\textbf{x}_{\mathbb{I}_{i,j,\kappa}},t), \quad \widehat{\mathbf{U}}_{R}=\mathbf{T}_{\mathbb{I}_{i,j}} \textbf{U}^n_{\mathbb{I}_{i,j,R}}(\textbf{x}_{\mathbb{I}_{i,j,\kappa}},t) .
    \end{equation}

    To close the spatial discretization, the remaining work reduces to evaluate the numerical flux $\mathbf{F}^x\left(\widehat{\mathbf{U}}_{L}, \widehat{\mathbf{U}}_{R}\right)$ by solving the simple one-dimensional Riemann problem 
    \begin{equation}
    \frac{\partial}{\partial t} \widehat{\mathbf{U}}+\frac{\partial}{\partial s} \widehat{\mathbf{F}}^x=\mathbf{0}, \quad \widehat{\mathbf{F}}^x=\mathbf{F}^x(\widehat{\mathbf{U}}), \quad \widehat{\mathbf{U}}(s, 0)=\left\{\begin{array}{ll}
    \widehat{\mathbf{U}}_{L}, & s<0 , \\
    \widehat{\mathbf{U}}_{R}, & s>0 .
    \end{array}\right.
    \end{equation}
    \subsection{HLL Riemman solver}
    Hereafter, the HLL \cite{harten1983upstream}\cite{toro1994restoration} approximate Riemann solver, which has been demonstrated to be robust and reliable for Euler equations, is briefly reviewed. The approximate Riemann solution for the discontinuous left and right state $\widehat{\mathbf{U}}_L$, $\widehat{\mathbf{U}}_R$ is
    \begin{equation}
    \tilde{\mathbf{U}}(x, t)=\left\{\begin{array}{clc}
    \widehat{\mathbf{U}}_{L}, & \text{if} & 0 \leq S_{L}, \\
    \widehat{\mathbf{U}}^{h l l}, & \text{if} & S_{L} \leq 0 , \\
    \widehat{\mathbf{U}}_{R}, & \text{if} & 0 \geq S_{R} ,
    \end{array}\right.
    \end{equation}
    and the corresponding interface flux is
    \begin{equation}
    \widehat{\mathbf{F}}^{x,hll}=\left\{\begin{array}{ccc}
    \widehat{\mathbf{F}}^x_{L}, & \text{if} & 0 \leq S_{L}, \\
    \frac{S_{R} \widehat{\mathbf{F}}^x_{L}-S_{L} \widehat{\mathbf{F}}^x_{R}+S_{L} S_{R}\left(\widehat{\mathbf{U}}_{R}-\widehat{\mathbf{U}}_{L}\right)}{S_{R}-S_{L}}, & \text { if } & S_{L} \leq 0 \leq S_{R}, \\
    \widehat{\mathbf{F}}^x_{R}, & \text{if} & 0 \geq S_{R}.
    \end{array}\right.
    \end{equation}
    where $\widehat{\mathbf{F}}^x_L = \mathbf{F}^x (\widehat{\mathbf{U}}_L)$ and $\widehat{\mathbf{F}}^x_R = \mathbf{F}^x (\widehat{\mathbf{U}}_R)$. The acoustic wave-speed $S_L$ and $S_R$ are evaluated by \cite{einfeldt1991godunov}\cite{batten1997choice}
    \begin{equation}
    \label{eq:wavespeed}
    {S_L} = \min ({u_L} - {c_L},\tilde u - \tilde c), \text{         } {S_R} = \max ({u_R} + {c_R},\tilde u + \tilde c),
    \end{equation}
    where the Roe-averaged quantities are defined as
    \begin{equation}
    \left\{ {\begin{array}{*{20}{c}}
    {{D_\rho } = \sqrt {\frac{\rho _R}{\rho _L}} } , \\
    {\tilde u = \frac{{u_L} + {u_R}{D_\rho }}{1 + {D_\rho }}} , \\
    {\tilde H = \frac{{H_L} + {H_R}{D_\rho }}{1 + {D_\rho }}} , \\
    {\tilde c = \sqrt {(\gamma  - 1)[\tilde H - \frac{1}{2}{{\tilde u}^2}]} } .
    \end{array}} \right.
    \end{equation}
\section{Numerical validations}

In this section, the proposed TENO schemes of third- to sixth-order are validated by solving the linear advection problem, the 1D and 2D Euler equations. If not mentioned otherwise, the third-order SSP Runge-Kutta scheme~\cite{gottlieb2001strong} is adopted for the temporal advancement with a constant CFL number 0.4. For all the simulations, the built-in parameters of TENO schemes are fixed and the results from CWENO and WENO schemes are provided for comparisons. {\color{black}All the present schemes are implemented in the research platform \cite{tsoutsanis2011weno}\cite{titarev2010weno}\cite{tsoutsanis2019stencil}\cite{tsoutsanis2021arbitrary}\cite{tsoutsanis2021cweno}, i.e., \texttt{UCNS3D} available at {\url{http://www.ucns3d.com}}, which has been extensively validated with both canonical and real-world simulations based on unstructured meshes.}

{\color{black}While HLL \cite{harten1983upstream} Riemman solver is a bit more dissipative than HLLC \cite{toro1994restoration} and Roe \cite{roe1981approximate}, it has lots
of merits, e.g., featuring better numerical robustness, free from Carbuncle problem and the
necessity of entropy fix. Since our main concern is the low-dissipation property offered by the new
TENO schemes, the choice of the Riemann solver does not affect our conclusion anyway and HLL provides a robust platform for these comparisons.}

\subsection{Accuracy order tests}

    We consider the 2D linear advection equation
    \begin{equation}
        \label{eq:accuracy_test_function}
        \frac{\partial u}{\partial t}+\frac{\partial u}{\partial x}+\frac{\partial u}{\partial y}=0,
    \end{equation}
    with the initial condition as
    \begin{equation}
        \label{eq:accuracy_test_function}
        u(x,y,0)=sin(2\pi x)\cdot sin(2\pi y).
    \end{equation}
    Periodic condition is imposed for all boundaries and the simulation is run until $t=1$.
    
    As shown in Table~\ref{Tab:accuracy_order_3},~\ref{Tab:accuracy_order_4},~\ref{Tab:accuracy_order_5}, and~\ref{Tab:accuracy_order_6}, the desired accuracy order is restored for all the considered TENO schemes with respect to both the $L^{\infty}$ and $L^{2}$ norms. Moreover, the numerical truncation errors from the nonlinear TENO schemes are identical to those from the linear schemes for all the resolutions indicating that the targeted linear schemes are restored exactly.

    \begin{table}[H]
        \centering
        \caption{Accuracy order test for third order schemes}
        \scriptsize
        \label{Tab:accuracy_order_3}
        \newcommand{\tabincell}[2]{\begin{tabular}{@{}#1@{}}#2\end{tabular}}
        \begin{tabular}{>{\centering\arraybackslash}m{1.2cm}
                        >{\centering\arraybackslash}m{1cm}
                        >{\centering\arraybackslash}m{1cm}
                        >{\centering\arraybackslash}m{1.4cm}
                        >{\centering\arraybackslash}m{1.4cm}
                        >{\centering\arraybackslash}m{1.4cm}
                        >{\centering\arraybackslash}m{1.4cm}}
        \hline
        & $h$ & $N_{e}\footnotemark$ & $L^{\infty}$ error & $L^{\infty}$ order\footnotemark & $L^{2}$ error & $L^{2}$ order \\ \hline
         Linear & 1/13  &  374   &  8.65E-02 &   -   & 4.58E-02 &   -   \\
                & 1/20  &  918   &  2.37E-02 &  2.88 & 1.30E-02 &  2.81 \\
                & 1/40  &  3674  &  3.33E-03 &  2.83 & 1.64E-03 &  2.98 \\
                & 1/80  &  14860 &  3.78E-04 &  3.12 & 2.01E-04 &  3.01 \\
                & 1/160 &  58330 &  5.14E-05 &  2.92 & 2.56E-05 &  3.01 \\
        \hline
         TENO3   & 1/13  &  374   &  8.65E-02 &   -   & 4.58E-02 &   -   \\
                & 1/20  &  918   &  2.37E-02 &  2.88 & 1.30E-02 &  2.81 \\
                & 1/40  &  3674  &  3.33E-03 &  2.83 & 1.64E-03 &  2.98 \\
                & 1/80  &  14860 &  3.78E-04 &  3.12 & 2.01E-04 &  3.01 \\
                & 1/160 &  58330 &  5.14E-05 &  2.92 & 2.56E-05 &  3.01 \\
        \hline
        \end{tabular}
    \end{table}
    \footnotetext[1]{$N_{e}$ the total number of elements}
    \footnotetext[2]{Both $L^{\infty}$ and $L^{2}$ order are calculated based on $N_{e}$}

    \begin{table}[H]
        \centering
        \caption{Accuracy order test for fourth order schemes}
        \scriptsize
        \label{Tab:accuracy_order_4}
        \newcommand{\tabincell}[2]{\begin{tabular}{@{}#1@{}}#2\end{tabular}}
        \begin{tabular}{>{\centering\arraybackslash}m{1.2cm}
                        >{\centering\arraybackslash}m{1cm}
                        >{\centering\arraybackslash}m{1cm}
                        >{\centering\arraybackslash}m{1.4cm}
                        >{\centering\arraybackslash}m{1.4cm}
                        >{\centering\arraybackslash}m{1.4cm}
                        >{\centering\arraybackslash}m{1.4cm}}
        \hline
        & $h$ & $N_{e}$ & $L^{\infty}$ error & $L^{\infty}$ order & $L^{2}$ error & $L^{2}$ order \\ \hline
         Linear & 1/13  &  374   &  7.99E-03 &   -   & 3.52E-03 &   -   \\
                & 1/20  &  918   &  2.26E-03 &  2.81 & 8.29E-04 &  3.22 \\
                & 1/40  &  3674  &  9.05E-05 &  4.64 & 3.28E-05 &  4.66 \\
                & 1/80  &  14860 &  6.59E-06 &  3.75 & 2.15E-06 &  3.90 \\
                & 1/160 &  58330 &  5.60E-07 &  3.61 & 1.57E-07 &  3.83 \\
        \hline
         TENO4   & 1/13  &  374   &  7.99E-03 &   -   & 3.52E-03 &   -   \\
                & 1/20  &  918   &  2.26E-03 &  2.81 & 8.29E-04 &  3.22 \\
                & 1/40  &  3674  &  9.05E-05 &  4.64 & 3.28E-05 &  4.66 \\
                & 1/80  &  14860 &  6.59E-06 &  3.75 & 2.15E-06 &  3.90 \\
                & 1/160 &  58330 &  5.60E-07 &  3.61 & 1.57E-07 &  3.83 \\
        \hline
        \end{tabular}
    \end{table}

    \begin{table}[H]
        \centering
        \caption{Accuracy order test for fifth order schemes}
        \scriptsize
        \label{Tab:accuracy_order_5}
        \newcommand{\tabincell}[2]{\begin{tabular}{@{}#1@{}}#2\end{tabular}}
        \begin{tabular}{>{\centering\arraybackslash}m{1.2cm}
                        >{\centering\arraybackslash}m{1cm}
                        >{\centering\arraybackslash}m{1cm}
                        >{\centering\arraybackslash}m{1.4cm}
                        >{\centering\arraybackslash}m{1.4cm}
                        >{\centering\arraybackslash}m{1.4cm}
                        >{\centering\arraybackslash}m{1.4cm}}
        \hline
        & $h$ & $N_{e}$ & $L^{\infty}$ error & $L^{\infty}$ order & $L^{2}$ error & $L^{2}$ order \\ \hline
         Linear & 1/13  &  374   &  9.16E-03 &   -   & 4.80E-03 &   -   \\
                & 1/20  &  918   &  1.06E-03 &  4.80 & 5.58E-04 &  4.79 \\
                & 1/40  &  3674  &  3.43E-05 &  4.95 & 1.59E-05 &  5.13 \\
                & 1/80  &  14860 &  9.68E-07 &  5.11 & 4.95E-07 &  4.97 \\
                & 1/160 &  58330 &  3.33E-08 &  4.93 & 1.61E-08 &  5.01 \\
        \hline
         TENO5   & 1/13  &  374   &  9.16E-03 &   -   & 4.80E-03 &   -   \\
                & 1/20  &  918   &  1.06E-03 &  4.80 & 5.58E-04 &  4.79 \\
                & 1/40  &  3674  &  3.43E-05 &  4.95 & 1.59E-05 &  5.13 \\
                & 1/80  &  14860 &  9.68E-07 &  5.11 & 4.95E-07 &  4.97 \\
                & 1/160 &  58330 &  3.33E-08 &  4.93 & 1.61E-08 &  5.01 \\
        \hline
        \end{tabular}
    \end{table}
    
    \begin{table}[H]
        \centering
        \caption{Accuracy order test for sixth order schemes}
        \scriptsize
        \label{Tab:accuracy_order_6}
        \newcommand{\tabincell}[2]{\begin{tabular}{@{}#1@{}}#2\end{tabular}}
        \begin{tabular}{>{\centering\arraybackslash}m{1.2cm}
                        >{\centering\arraybackslash}m{1cm}
                        >{\centering\arraybackslash}m{1cm}
                        >{\centering\arraybackslash}m{1.4cm}
                        >{\centering\arraybackslash}m{1.4cm}
                        >{\centering\arraybackslash}m{1.4cm}
                        >{\centering\arraybackslash}m{1.4cm}}
        \hline
        & $h$ & $N_{e}$ & $L^{\infty}$ error & $L^{\infty}$ order & $L^{2}$ error & $L^{2}$ order \\ \hline
         Linear & 1/13  &  374   &  7.94E-04 &   -   & 3.17E-04 &   -   \\
                & 1/20  &  918   &  8.25E-05 &  5.04 & 2.60E-05 &  5.57 \\
                & 1/40  &  3674  &  1.21E-06 &  6.08 & 3.50E-07 &  6.22 \\
                & 1/80  &  14860 &  2.21E-08 &  5.73 & 5.48E-09 &  5.95 \\
                & 1/160 &  58330 &  3.72E-10 &  5.97 & 9.69E-11 &  5.90 \\
        \hline
         TENO6   & 1/13  &  374   &  7.94E-04 &   -   & 3.17E-04 &   -   \\
                & 1/20  &  918   &  8.25E-05 &  5.04 & 2.60E-05 &  5.57 \\
                & 1/40  &  3674  &  1.21E-06 &  6.08 & 3.50E-07 &  6.22 \\
                & 1/80  &  14860 &  2.21E-08 &  5.73 & 5.48E-09 &  5.95 \\
                & 1/160 &  58330 &  3.72E-10 &  5.97 & 9.69E-11 &  5.90 \\
        \hline
        \end{tabular}
    \end{table}

\subsection{{\color{black}2D isentropic vortex evolution}}

{\color{black}We consider the 2D isentropic vortex evolution problem in the computational domain of $[0,10]\times[0,10]$ with periodic boundary conditions on all sides \cite{hu1999weighted}. 
The unperturbed initial condition is $\left(\rho,u,v,p\right)=\left(1,1,1,1\right)$, and the vortex perturbations (no pertubation in entropy $S=p/{\rho}^{\gamma}$) are given as
\begin{equation}\begin{split}
& \left(\delta u,\delta v\right)=\frac{\epsilon}{2\pi}e^{0.5\left(1-r^{2}\right)}\left(-\left(y-5\right),\left(x-5\right)\right),\\
&\delta T=-\frac{\left(\gamma-1\right)\epsilon^{2}}{8\gamma\pi^{2}} e^{\left(1-r^{2}\right)},
\end{split}
\end{equation}
where the temperature is defined as $T=p/\rho$, the vortex strength is $\epsilon=5$, the gas constant $\gamma=1.4$ and $r^2=(x-5)^2+(y-5)^2$. The simulation time is $t=10$.}

    \begin{figure}[H]
        \centering
        \includegraphics[width=0.6\textwidth,trim=4 4 4 4,clip]{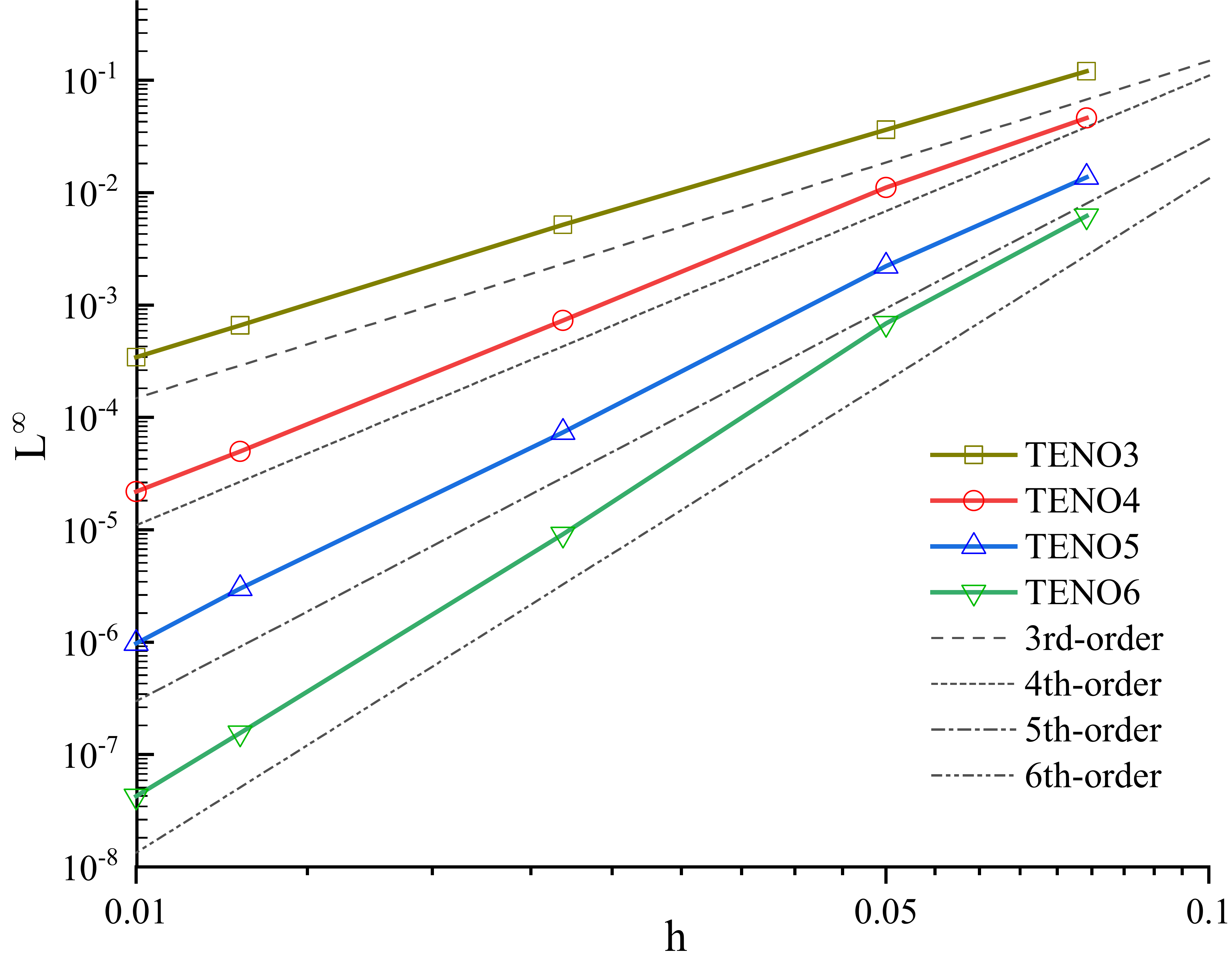}
        \caption{{\color{black}2D isentropic vortex evolution: convergence statistics of the $L^{\infty}$ error from all the considered TENO schemes.}}
        \label{Fig_order}
    \end{figure}

{\color{black}As shown in Fig.~\ref{Fig_order}, the desired convergence order is achieved for all the considered TENO schemes.}
\subsection{2D solid body rotation problem}

    Next we consider the solid body rotation problem from Leveque \cite{leveque1996high}. The computational domain is $[0,1]\times[0,1]$ with periodic boundary conditions applied for all directions. Three bodies are considered in this case, i.e., a cosine function (SB01), a sharp cone (SB02) and a slotted cylinder (SB03). The functions, center locations and radius to describe the three bodies are provided in Table~\ref{Tab:2d_rotation}. The bodies are advected by an initial velocity $\textbf{v}(x,y)=(0.5-y,x-0.5)$ around the center of the domain, i.e. $(x_c,y_c) = (0.5,0.5)$. At $t=2\pi$, the analytical profile coincides with the initial solution. A uniform triangular mesh of edge length $h=1/80$ is used and the total number of mesh elements is 14860.

    \begin{table}[H]
        \renewcommand{\arraystretch}{1.7}
        \centering
        \caption{Initial condition for the 2D solid body rotation problem}
        \scriptsize
        \label{Tab:2d_rotation}
        \newcommand{\tabincell}[2]{\begin{tabular}{@{}#1@{}}#2\end{tabular}}
        \begin{tabular}{>{\centering\arraybackslash}m{1.2cm}
                        >{\centering\arraybackslash}m{6cm}
                        >{\centering\arraybackslash}m{2cm}
                        >{\centering\arraybackslash}m{1.5cm}} 
        \hline
        & $\textbf{f}(x,y)$ & center $(x_0,y_0)$  & radius $r_0$\\
        \hline
         SB01 & $\dfrac{1}{4}(1+\cos(\pi r(x,y)))$  &  $(0.25,0.5)$   & 0.15\\
         SB02 & $1-(1/r_0)\sqrt{(x-x_0)^2+(y-y_0)^2}$ &  $(0.5,0.25)$  & 0.15 \\
         SB03 &$\left\{ {\begin{array}{*{2}{c}}{1,}&{\text{if }|x-x_0|\ge 0.025\text{ or }y\ge 0.85,} \\{0,}&{\text{otherwise.}}\end{array}} \right. $ &  $(0.5,0.75)$  & 0.15 \\
        \hline
        \end{tabular}
    \end{table}
    
\begin{figure}[H]
        \centering
        \includegraphics[width=1.0\textwidth,trim=4 4 4 4,clip]{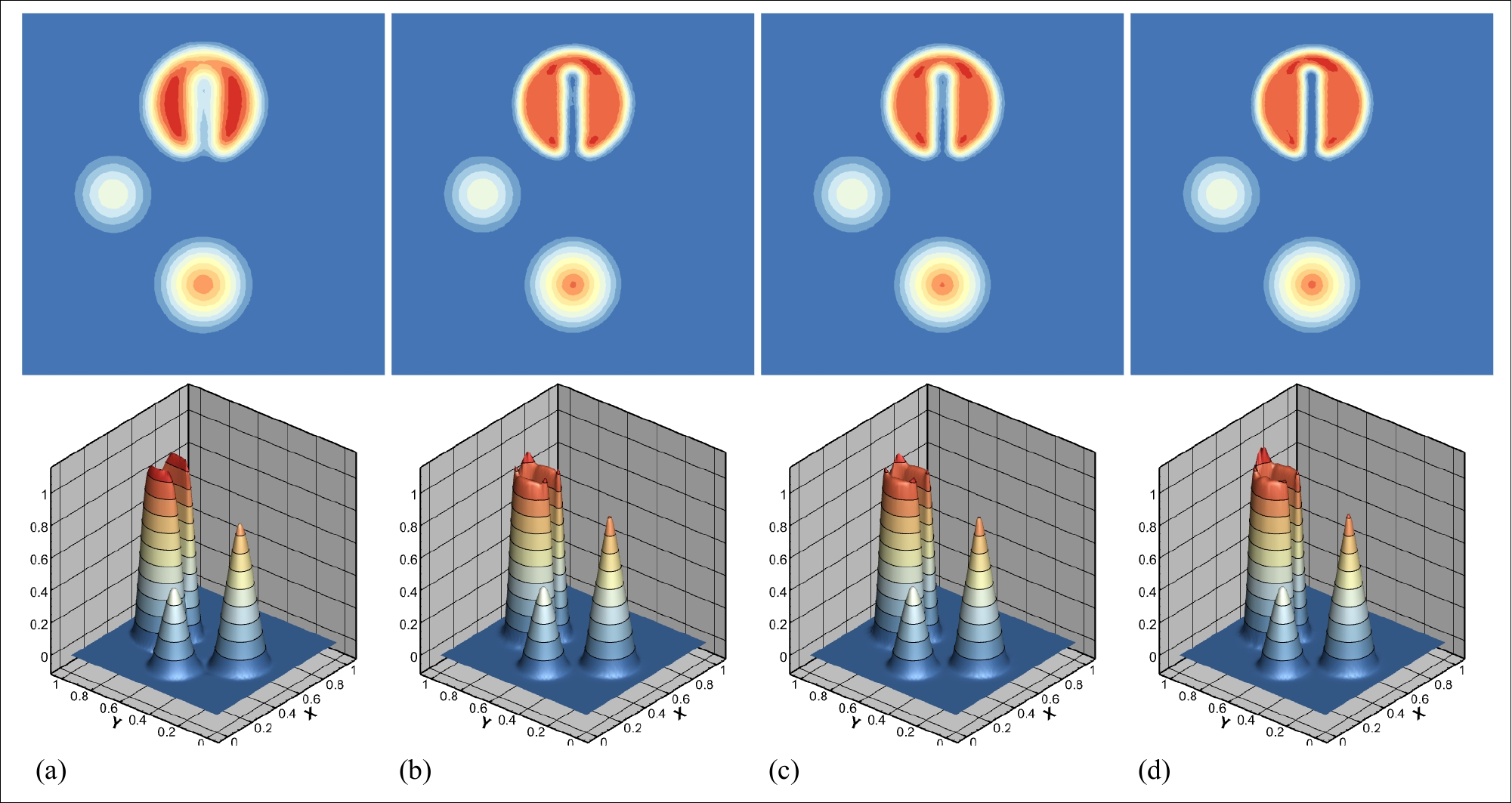}
        \caption{2D solid body rotation problem: results from the (a) TENO3, (b) TENO4, (c) TENO5 and (d) TENO6 schemes after one rotation evolution at $t=2\pi$. The uniform triangular mesh of edge length $h=1/80$ is used.}
        \label{Fig:2d_rotation_01}
    \end{figure}

\begin{figure}[H]
        \centering
        \includegraphics[width=1.0\textwidth,trim=4 4 4 4,clip]{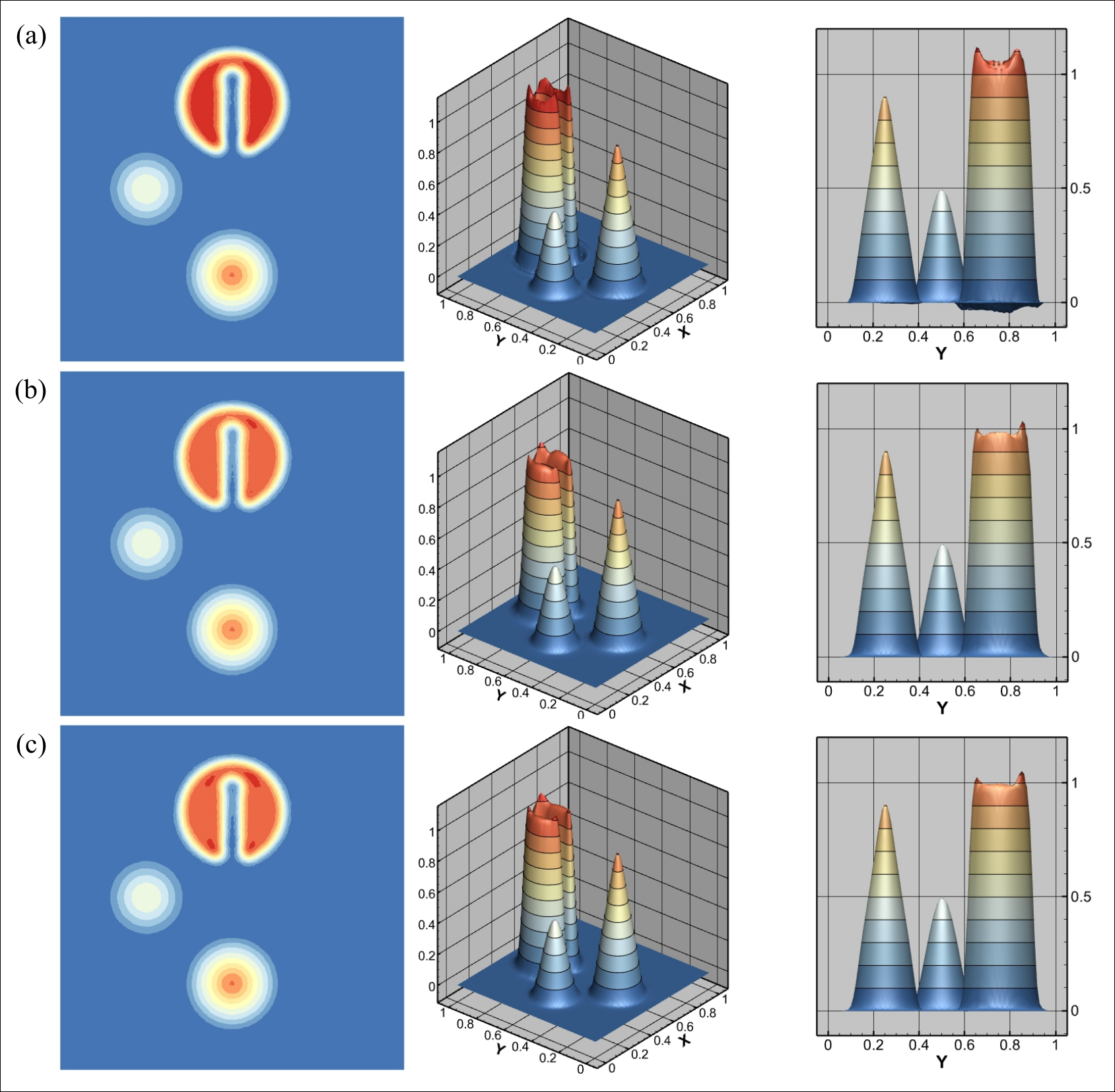}
        \caption{2D solid body rotation problem: results from the (a) WENO5, (b) CWENO5 and (c) TENO5 schemes after one rotation evolution at $t=2\pi$. The uniform triangular mesh of edge length $h=1/80$ is used.}
        \label{Fig:2d_rotation_02}
    \end{figure}

    As shown in Fig.~\ref{Fig:2d_rotation_01}, the high-order TENO schemes, i.e. TENO4, TENO5 and TENO6, show much less numerical dissipation than the third-order TENO3 in capturing the slotted cylinder. Moreover, as demonstrated in Fig.~\ref{Fig:2d_rotation_02}, TENO5 generates less numerical overshoots than WENO5 and less dissipation than CWENO5 in preserving the shape of the cylinder.

\subsection{1D shock-tube problems}
    The initial state for the 1D shock tube (ST) problems \cite{dumbser2017central} is
    \begin{equation}
        \label{eq:shock_tube_initial states}
            (\rho ,u,p) = \left\{ {\begin{array}{*{20}{c}}
            {(\rho_L,u_L,p_L),}&{\text{if }0 \le x < x_d ,} \\
            {(\rho_R,u_R,p_R),}&{\text{if }x_d \le x \le 1 ,}
        \end{array}} \right. 
    \end{equation}
    where $L$ and $R$ denotes the left and right states regarding to the location of discontinuity at $x = x_{d}$. The final simulation time is $t = t_{end}$. We consider three variations of this problem, i.e. ST01/02/03, and all the parameters are given in Table~\ref{Tab:sod_initial_states}. The computational domain is $[0,0.2]\times[0,1]$. For this case, the uniform Cartesian mesh is employed with an effective resolution of $h=1/100$ and the conservative variables are deployed for the high-order reconstruction due to the simplicity. 

    \begin{table}[H]
        \centering
        \caption{Parameters of the 1D shock-tube problems. }
        \scriptsize
        \label{Tab:sod_initial_states}
        \newcommand{\tabincell}[2]{\begin{tabular}{@{}#1@{}}#2\end{tabular}}
        \begin{tabular}{>{\centering\arraybackslash}m{1.5cm}
                        >{\centering\arraybackslash}m{1cm}
                        >{\centering\arraybackslash}m{1cm}
                        >{\centering\arraybackslash}m{1cm}
                        >{\centering\arraybackslash}m{1cm}
                        >{\centering\arraybackslash}m{1cm}
                        >{\centering\arraybackslash}m{1.1cm}                    >{\centering\arraybackslash}m{1cm}
                        >{\centering\arraybackslash}m{1cm}}
        \hline
        & $\rho_L$ & $u_{L}$ & $p_{L}$ & $\rho_{R}$ & $u_{R}$ & $p_{R}$ & $t_{end}$ & $x_{d}$ \\ \hline
         ST01 (Sod)  & 1.0     &  0.0     &  1.0     & 1.0     & 0.125   &  0.0     & 0.2   &  0.5   \\
         ST02 (Lax)  & 0.445   &  0.698   &  3.528   & 3.528   & 0.5     &  0.0     & 0.14  &  0.5   \\
         ST03        & 1.0     &  0.0     &  1000    & 1000    & 1.0     &  0.0     & 0.012 &  0.6   \\
        \hline
        \end{tabular}
    \end{table}

\begin{figure}[H]
        \centering
        \includegraphics[width=1.0\textwidth]{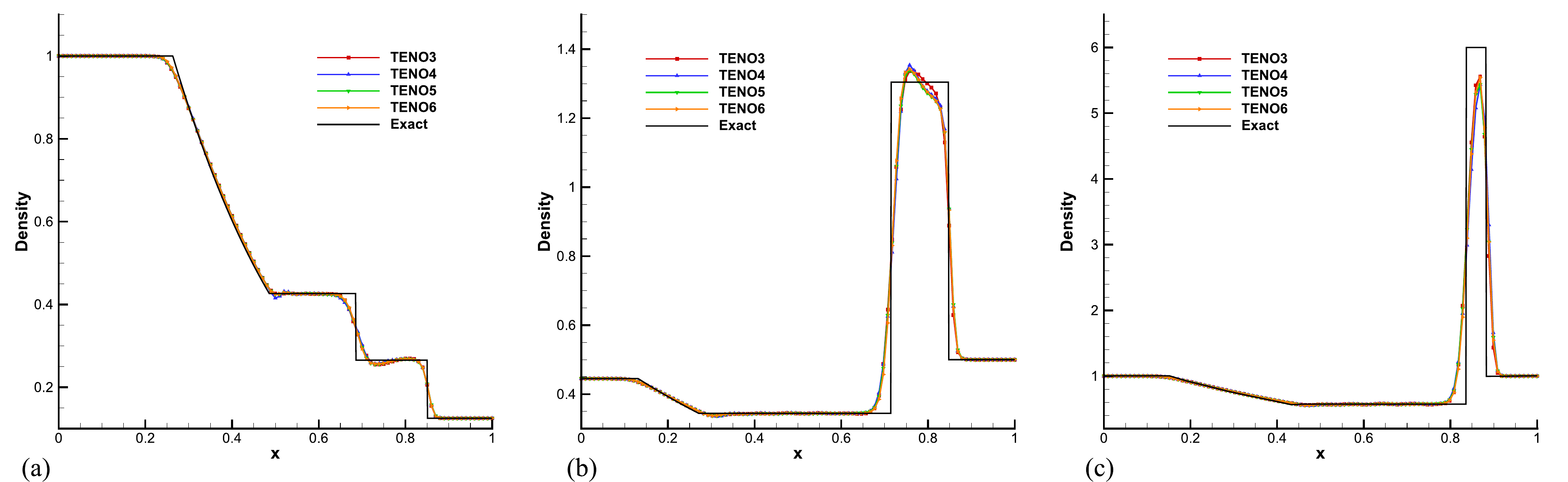}
        \caption{1D shock-tube problems: results of (a) ST01, (b) ST02 and (c) ST03 with TENO schemes of different orders. The mesh resolution is $h=1/100$.}
        \label{Fig:1d_shock_tube}
    \end{figure}

The computed density profiles are given in Fig.~\ref{Fig:1d_shock_tube}. Overall speaking, the discontinuities are captured without spurious oscillations by the TENO schemes of various orders. There are slight overshoots near the contact waves for the  ST02 problem, which are also observed even for the results from other finite-difference low-dissipation schemes on structured meshes \cite{fu2016family}.

\begin{figure}[H]
        \centering
        \includegraphics[width=1.0\textwidth]{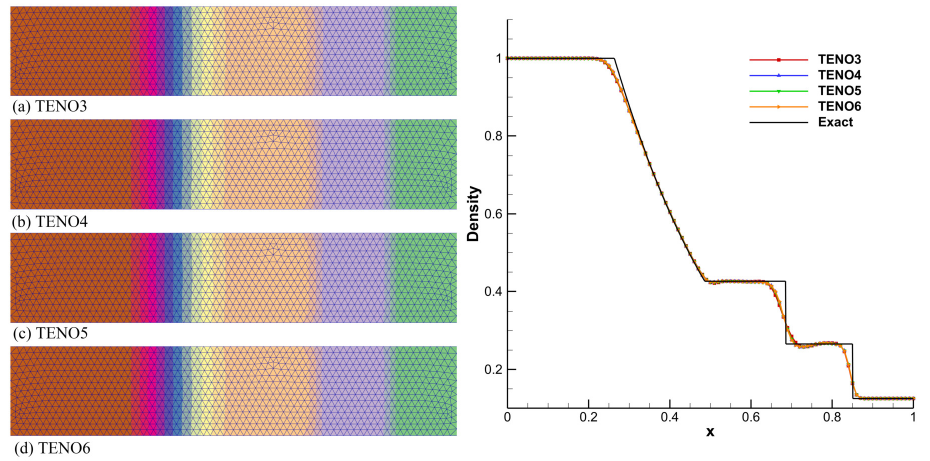}
        \caption{{\color{black}1D shock-tube problem (ST01): solutions from all the considered TENO schemes. Left: density contours and mesh typologies; right: density profiles.}}
        \label{Fig:sod2}
    \end{figure}

\begin{figure}[H]
        \centering
        \includegraphics[width=1.0\textwidth]{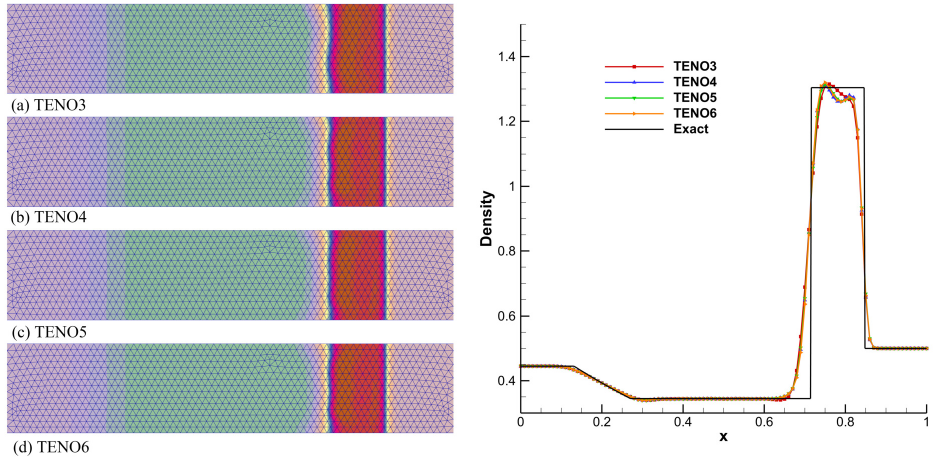}
       \caption{ {\color{black}1D shock-tube problem (ST02): solutions from all the considered TENO schemes. Left: density contours and mesh typologies; right: density profiles.}}
        \label{Fig:lax2}
    \end{figure}

{\color{black}In order to verify the capability of the present TENO schemes in maintaining the 1D symmetry, both the ST01 and ST02 problems are simulated additionally on a genuinely unstructured mesh with 1967 uniform triangular elements. As shown in Fig.~\ref{Fig:sod2} and Fig.~\ref{Fig:lax2}, the present TENO schemes with various orders preserve the 1D symmetry well.}

\subsection{1D shock-density wave interaction}

The initial condition is given as \cite{Shu1989}
    \begin{equation}
        \label{eq:shuosher}
        (\rho ,u,p) = \left\{ {\begin{array}{*{20}{c}}
        {(3.857,2.629,10.333) , }&{\text{if } 0\le x < 1 ,}\\
        {(1 + 0.2\sin (5(x-5)),0,1) , }&{\text{if } 1\le x \le 10 .}
        \end{array}} \right.
    \end{equation}
    The computational domain is $[0,2]\times[0,10]$ with $18,282$ uniformly distributed triangular mesh cells, i.e. $h\approx1/20$, and the final evolution time is $t = 1.8$. The ``exact" solution is obtained by the 1D fifth-order finite-difference WENO5-JS scheme with $2000$ cells.
      
\begin{figure}[H]
        \centering
        \includegraphics[width=1.0\textwidth,trim=2 2 2 2,clip]{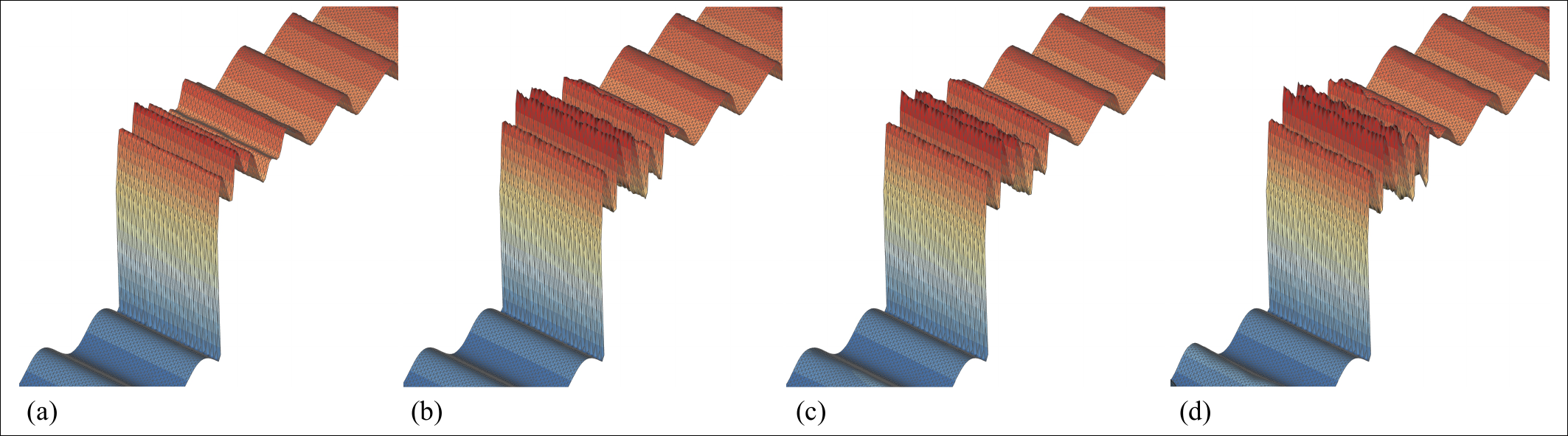}
        \caption{1D shock-density wave interaction: density distributions from the (a) TENO3, (b) TENO4, (c) TENO5 and (d) TENO6 schemes. Also shown is the mesh topology on top of the density distributions. The uniformly distributed triangular mesh with $h\approx1/20$ is adopted.}
        \label{Fig:2d_shu_01}
    \end{figure}

\begin{figure}[H]
        \centering
        \includegraphics[width=1.0\textwidth,trim=2 2 2 2,clip]{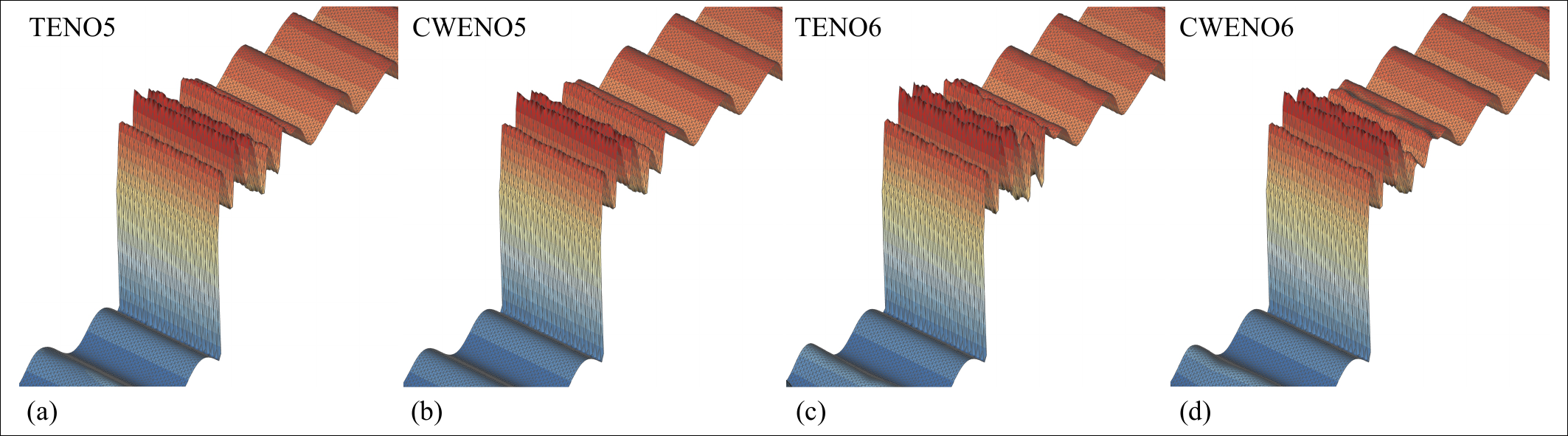}
        \caption{1D shock-density wave interaction: density distributions from the (a) TENO5, (b) CWENO5, (c) TENO6 and (d) CWENO6 schemes. Also shown is the mesh topology on top of the density distributions. The uniformly distributed triangular mesh with $h\approx1/20$ is adopted.}
        \label{Fig:2d_shu_02}
    \end{figure}

\begin{figure}[H]
        \centering
        \includegraphics[width=1.0\textwidth]{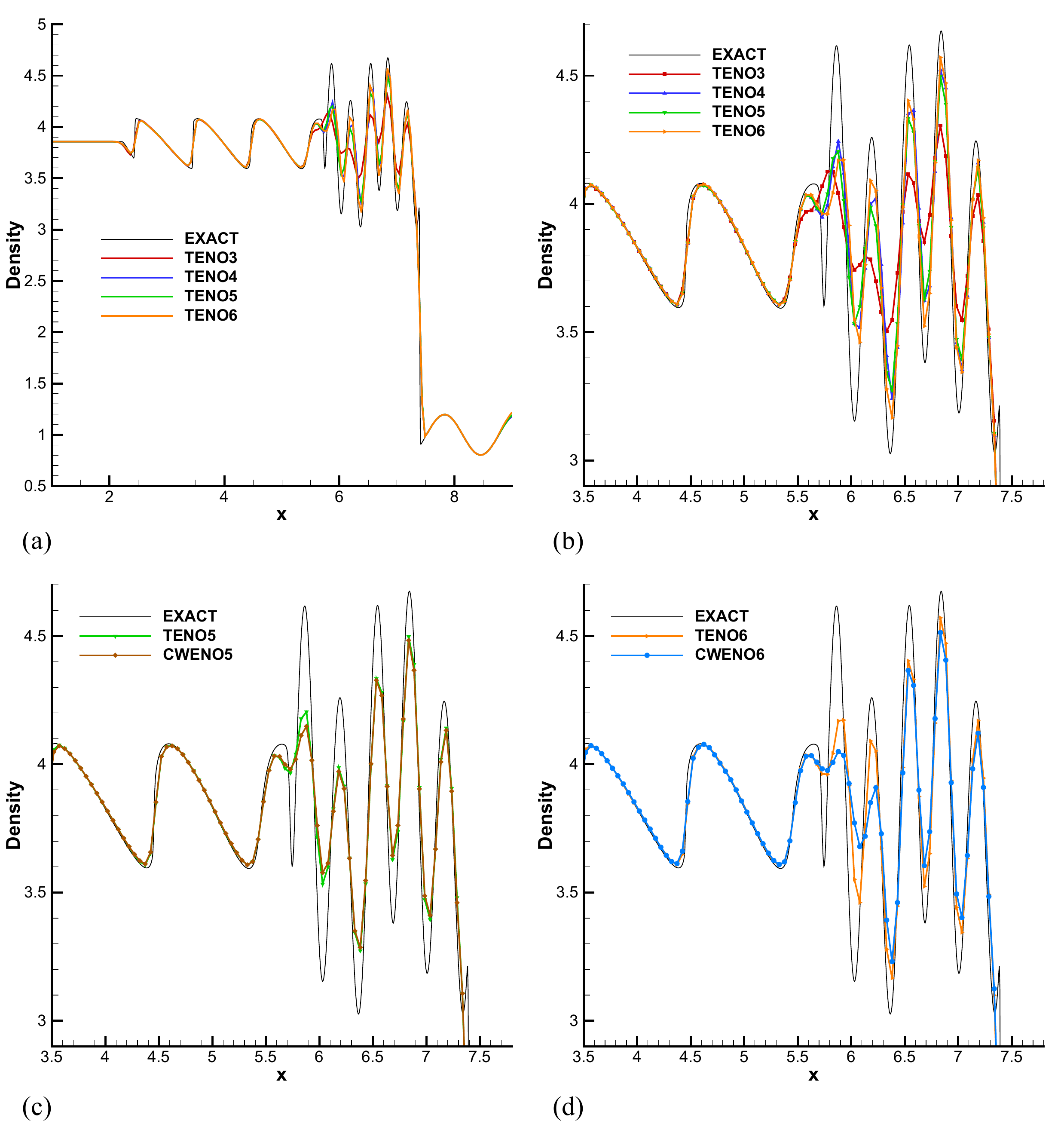}
        \caption{1D shock-density wave interaction: (a) density and (b) zoomed-in view of density profiles from TENO schemes of different orders. Also included are the comparisons between (c) TENO5 and CWENO5, and (d) TENO6 and CWENO6 for density profiles. The uniformly distributed triangular mesh with $h\approx1/20$ is adopted.}
        \label{Fig:2d_shu_03}
    \end{figure}

    As shown in Fig.~\ref{Fig:2d_shu_01},~\ref{Fig:2d_shu_02}, and~\ref{Fig:2d_shu_03}, with increasing accuracy order, the TENO schemes resolve the high-wavenumber fluctuations significantly better while maintaining the sharp shock-capturing capability. No artificial overshoots and oscillations are observed. In particular, the TENO5 and TENO6 schemes perform better than the CWENO schemes of the same accuracy order, suggesting less numerical dissipation.

\begin{figure}[H]
        \centering
        \includegraphics[width=1.0\textwidth]{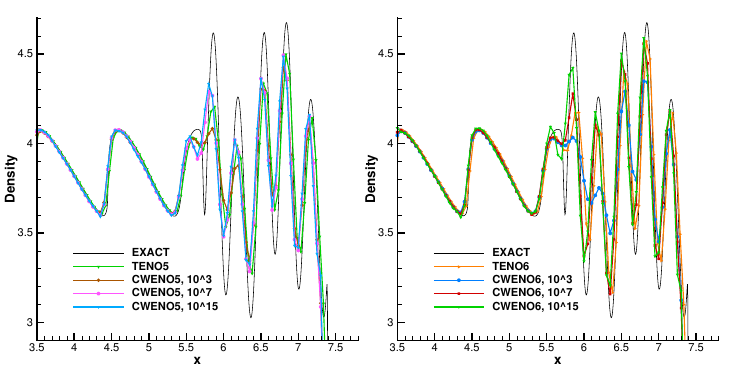}
        \caption{{\color{black}1D shock-density wave interaction: zoomed-in views of the density profiles from the CWENO schemes with different $d_K^{'}$ values and the TENO schemes. The uniformly distributed triangular mesh with $h\approx1/20$ is adopted.}}
        \label{Fig:shu2}
    \end{figure}

{\color{black}
As reported by Tsoutsanis and Dumbser \cite{tsoutsanis2021arbitrary}, different $d_K^{'}$ values can be adopted for CWENO schemes. In order to comprehensively compare the performance of TENO shcmes with that of CWENO schemes taking different values of $d_K^{'}$, Fig.~\ref{Fig:shu2} shows the density profiles of the present TENO schemes and the CWENO schemes with $d_K^{'} = 10^3, 10^7$, and $10^{15}$. 
It can be seen that although the CWENO scheme with larger $d_K^{'}$ performs better in resolving the high-wavenumber fluctuations, the robustness is sacrificed, as evidenced by additional numerical experiments. Especially when handling the problems involving strong discontinuities, such as the double Mach reflection problem, the simulation typically fails with too large $d_K^{'}$. As a compromise of performance and robustness, we adopt the uniform optimal $d_K^{'}$ value, i.e., $10^4$ for CWENO schemes, in the present study if not mentioned otherwise.
}

\subsection{2D explosion problems}
    Two 2D explosion problems, referred as EP01 and EP02, are considered. The computational domain is given by a unit circle of radius R = 1 centered at $[1,1]$. The computational domain is discretized by 18,464 triangular meshes, i.e. with an effective resolution of $h\approx1/50$. The initial condition for problem EP01 is
    \begin{equation}
        \label{eq:bw_ep1}
        (\rho,u,v,p)(\textbf{r},0)=
        \begin{cases}
        (1.0,0,0,1.0), & \text{if $ \Vert \textbf{r} \Vert \leqslant 0.5$},\\
        (0.125,0,0,0.1), & \text{otherwise}.
        \end{cases}
    \end{equation}
    The final simulation time is $t = 0.2$.

\begin{figure}[H]
        \centering
        \includegraphics[width=1.0\textwidth,trim=2 2 2 2,clip]{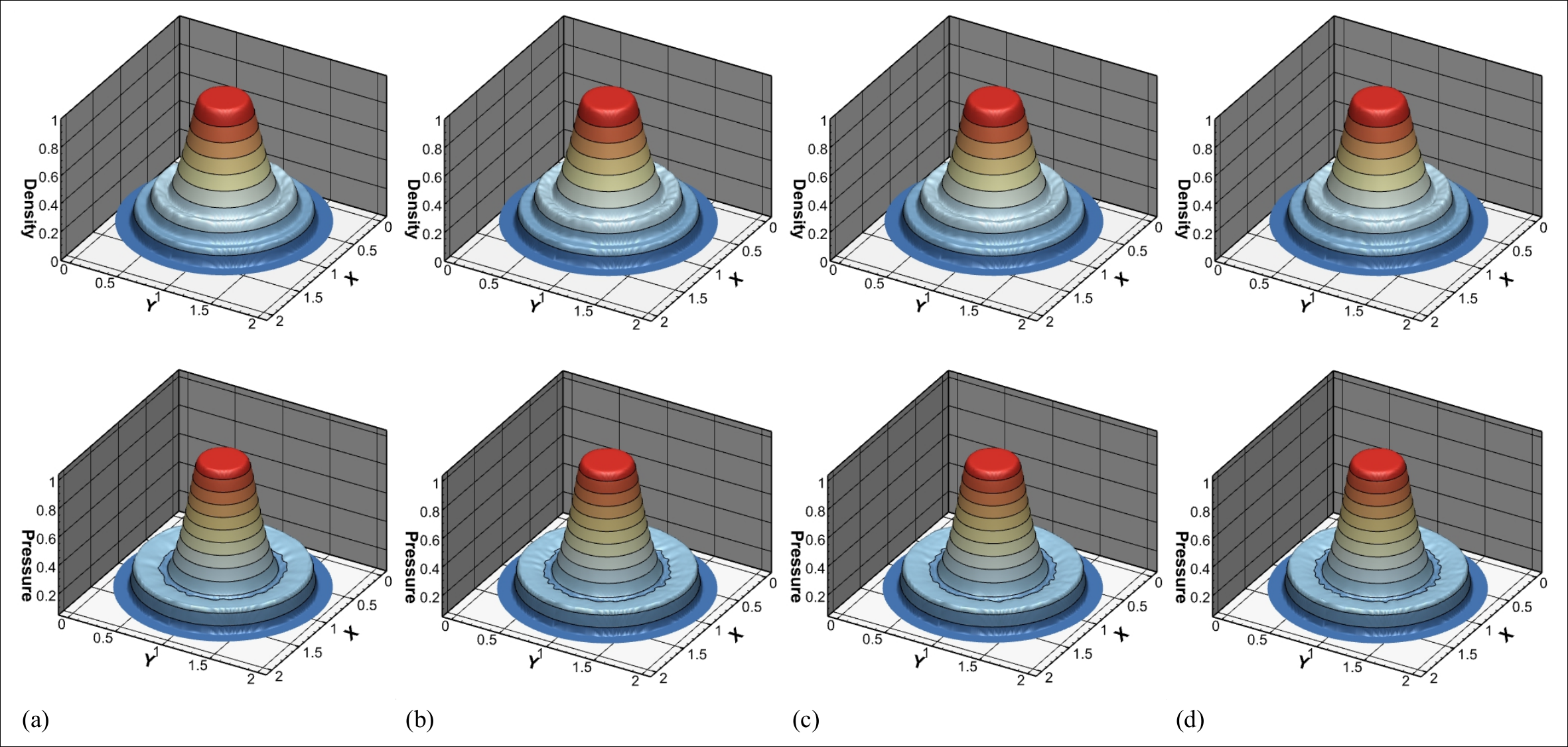}
        \caption{EP01: results from the (a) TENO3, (b) TENO4, (c) TENO5 and (d) TENO6 schemes. The upper row shows the density profiles and the bottom row shows the pressure profiles. A uniform triangular mesh with edge length $h\approx1/50$ is used.}
        \label{Fig:2d_explosion_ep01_01}
    \end{figure}

\begin{figure}[H]
        \centering
        \includegraphics[width=1.0\textwidth,trim=2 2 2 2,clip]{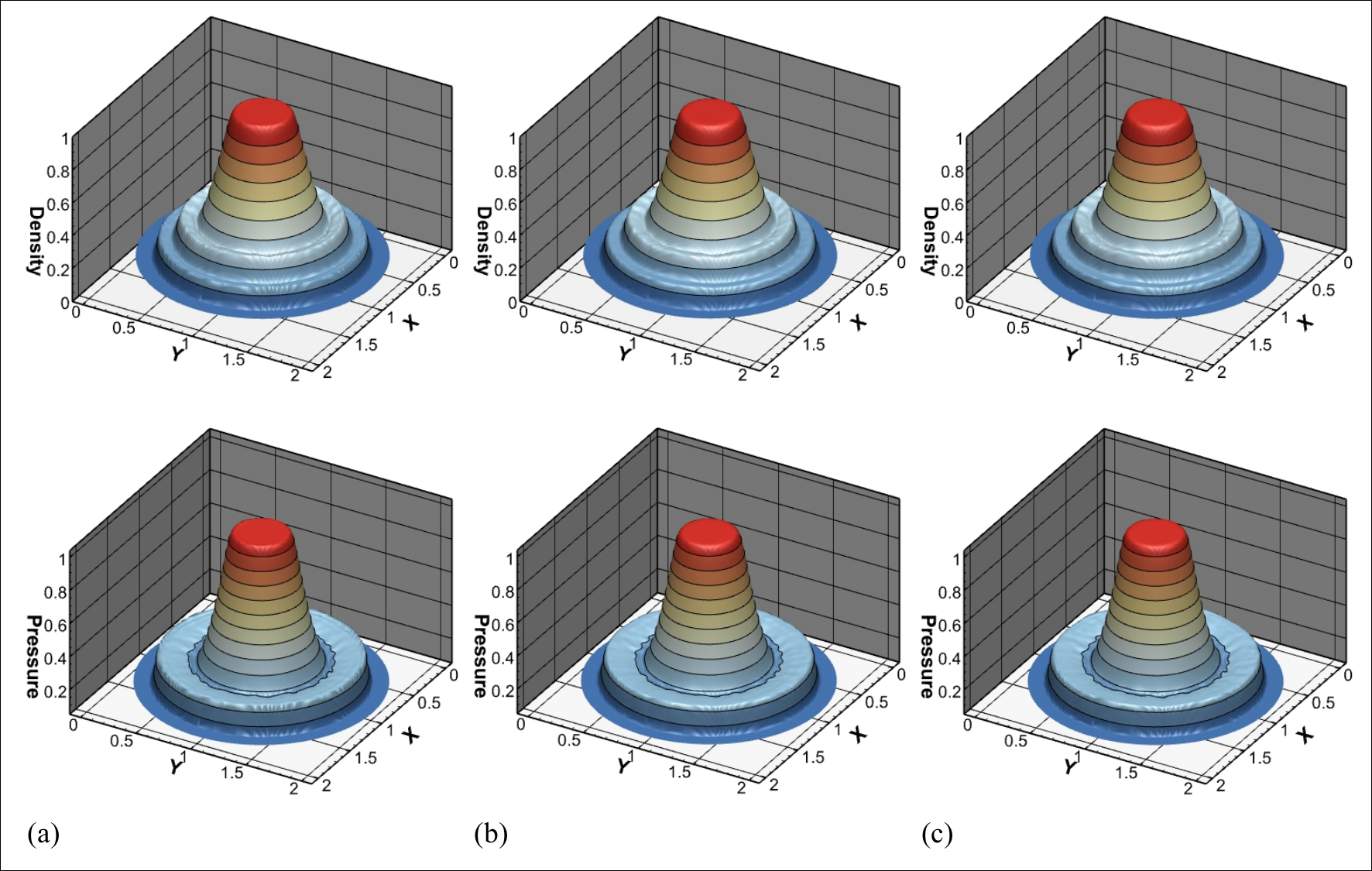}
        \caption{EP01: results from the (a) WENO5, (b) CWENO5 and (c) TENO5 schemes. The upper row shows the density profiles and the bottom row shows the pressure profiles. A uniform triangular mesh with edge length $h\approx1/50$ is used.}
        \label{Fig:2d_explosion_ep01_02}
    \end{figure}

\begin{figure}[H]
        \centering
        \includegraphics[width=1.0\textwidth]{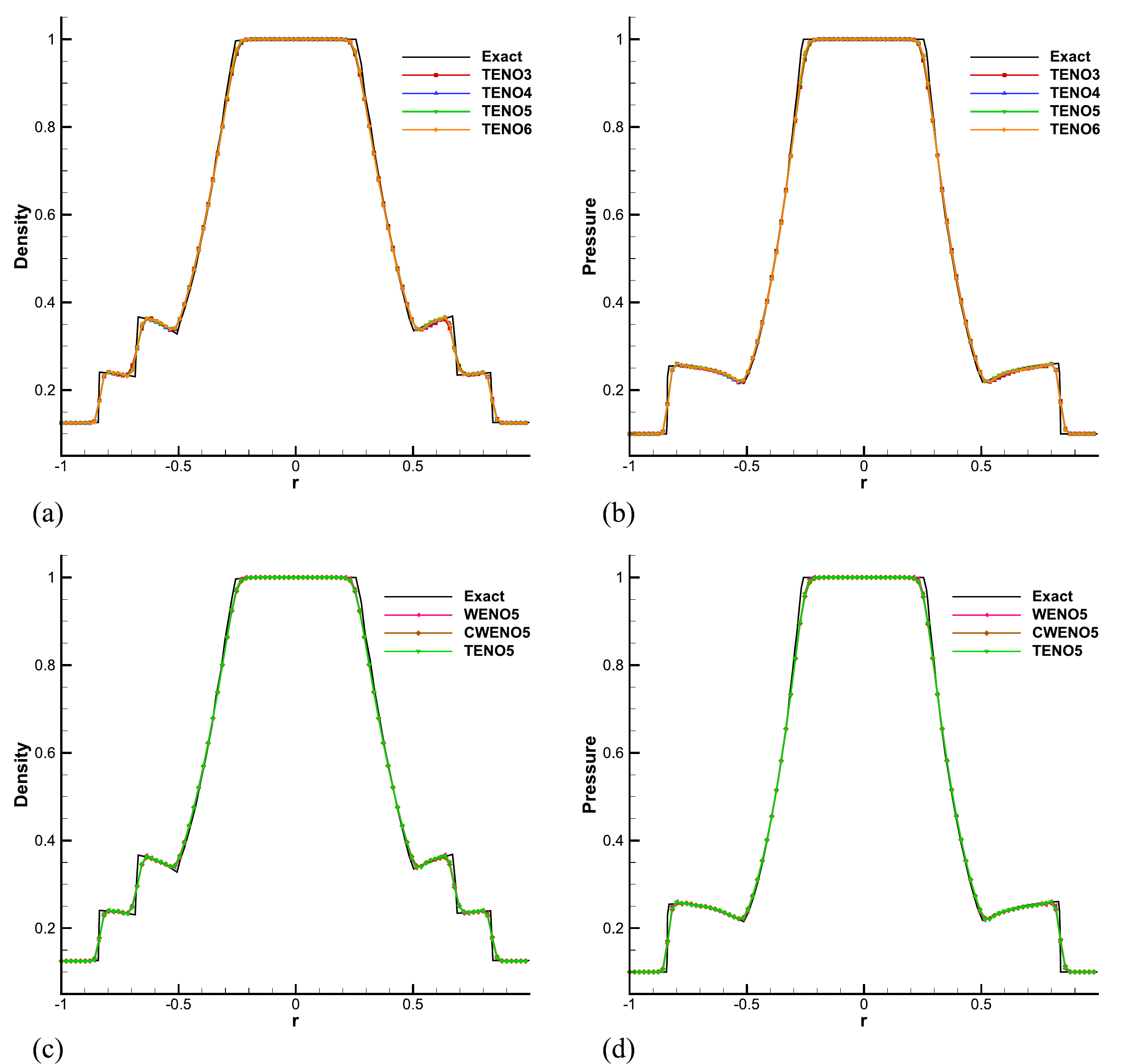}
        \caption{EP01: (a) density and (b) pressure profiles from TENO schemes of different orders. (c) density and (d) pressure  profiles from the TENO5, WENO5 and CWENO5 schemes. A uniform triangular mesh with edge length $h\approx1/50$ is used.}
        \label{Fig:2d_explosion_ep01_03}
    \end{figure}

    As shown in Fig.~\ref{Fig:2d_explosion_ep01_01},~\ref{Fig:2d_explosion_ep01_02}, and~\ref{Fig:2d_explosion_ep01_03}, for both the density and pressure profiles of EP01, the present TENO schemes with various orders preserve the monotonicity pretty well as well as the corresponding WENO5 and CWENO5 schemes.


    The initial condition for problem EP02 is
    \begin{equation}
        \label{eq:bw_ep2}
        (\rho,u,v,p)(\textbf{r},0)=
        \begin{cases}
        (1.0,0,0,2.0), & \text{if $ \Vert \textbf{r} \Vert \leqslant 0.5$},\\
        (1.0,0,0,1.0), & \text{otherwise}.
        \end{cases}
    \end{equation}
    The final simulation time is $t = 0.2 $.

\begin{figure}[H]
        \centering
        \includegraphics[width=1.0\textwidth,trim=2 2 2 2,clip]{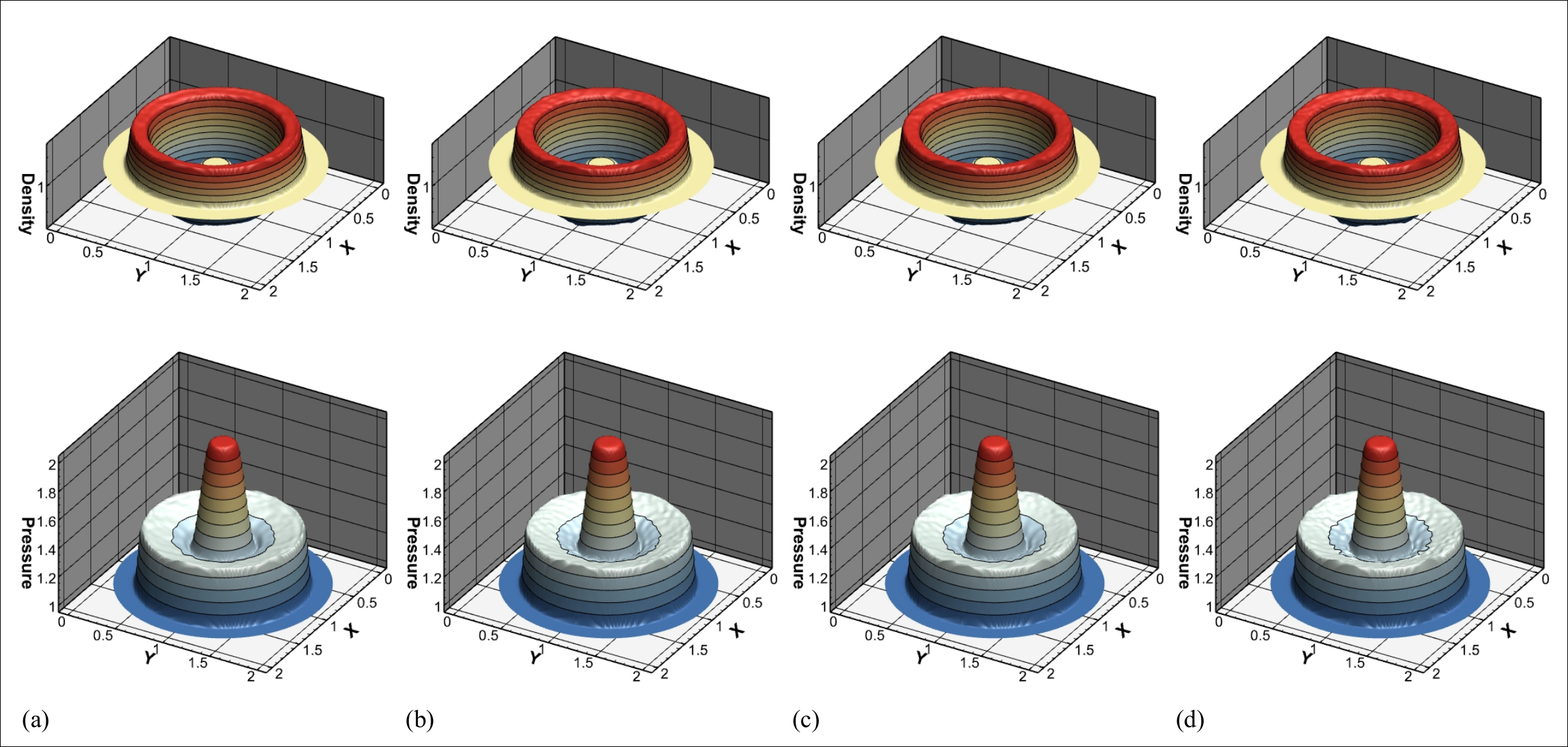}
        \caption{EP02: results from the (a) TENO3, (b) TENO4, (c) TENO5 and (d) TENO6 schemes. The upper row shows the density profiles and the bottom row shows the pressure profiles. A uniform triangular mesh with edge length $h\approx1/50$ is used.}
        \label{Fig:2d_explosion_ep02_01}
    \end{figure}

\begin{figure}[H]
        \centering
        \includegraphics[width=1.0\textwidth,trim=2 2 2 2,clip]{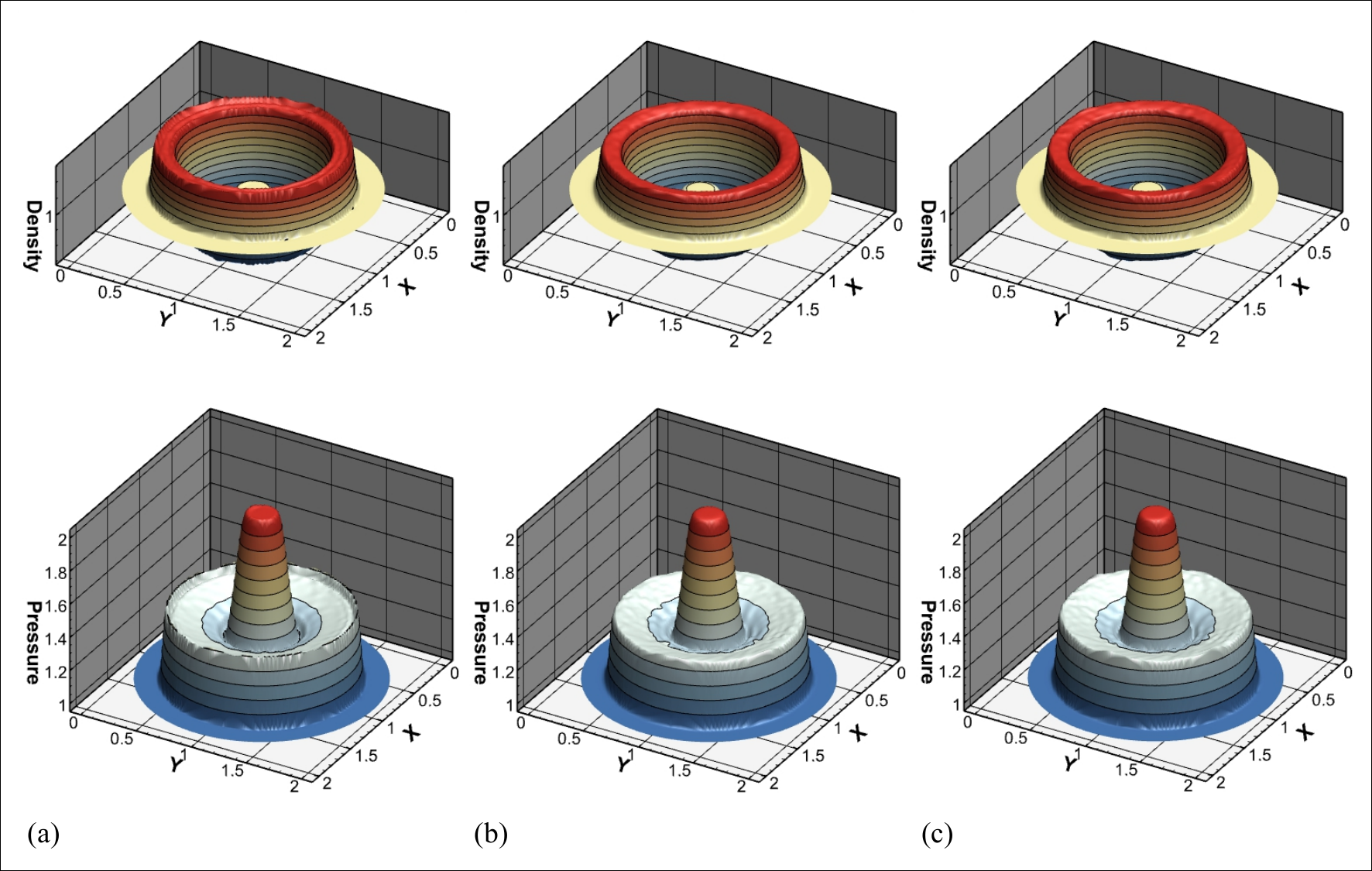}
        \caption{EP02: results from the (a) WENO5, (b) CWENO5 and (c) TENO5 schemes. The upper row shows the density profiles and the bottom row shows the pressure profiles. A uniform triangular mesh with edge length $h\approx1/50$ is used.}
        \label{Fig:2d_explosion_ep02_02}
    \end{figure}

\begin{figure}[H]
        \centering
        \includegraphics[width=1.0\textwidth]{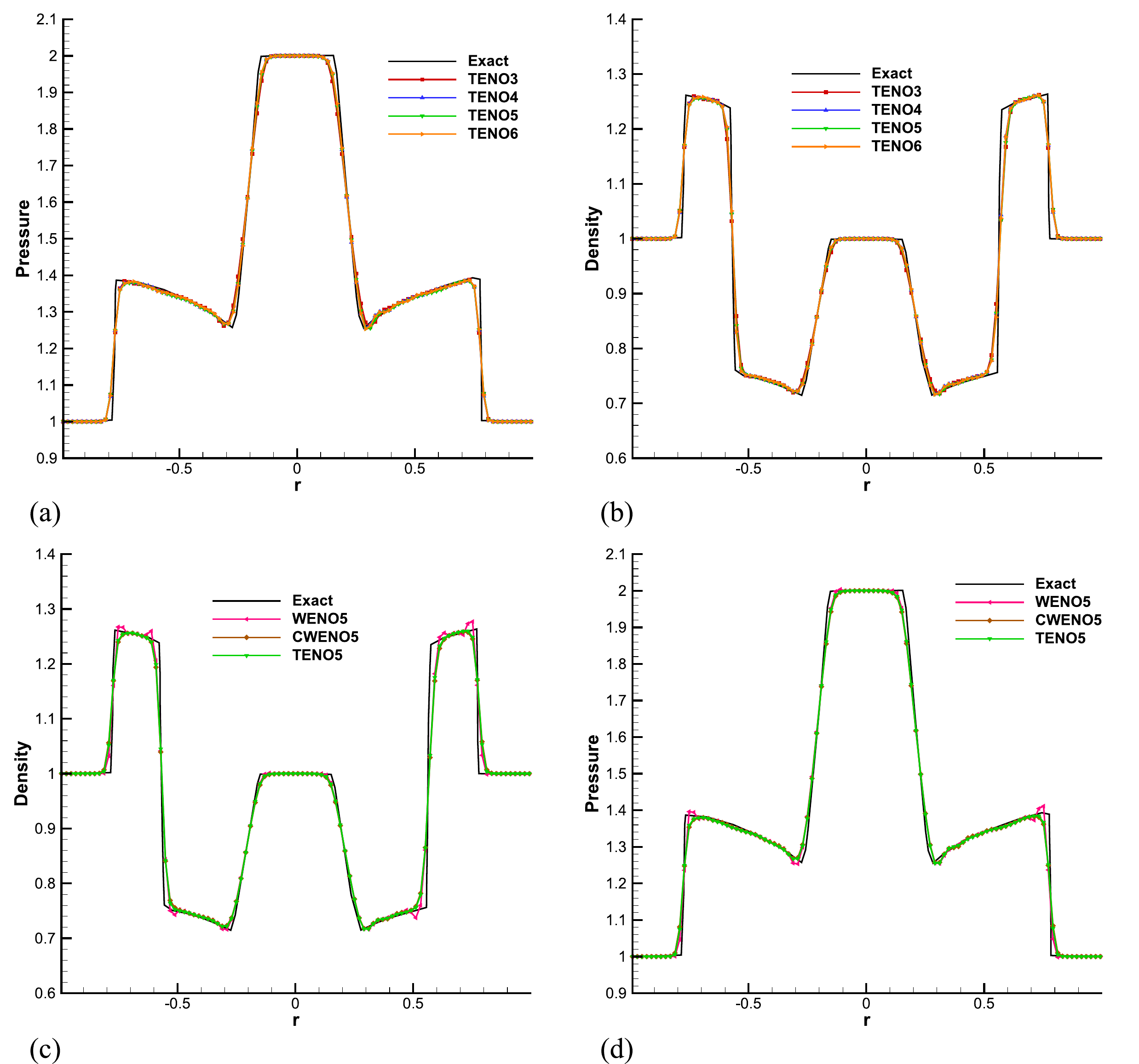}
        \caption{EP02: (a) density and (b) pressure profiles from TENO schemes of different orders; (c) density and (d) pressure profiles from the TENO5, WENO5 and CWENO5 schemes. A uniform triangular mesh with edge length $h\approx1/50$ is used.}
        \label{Fig:2d_explosion_ep02_03}
    \end{figure}

    In terms of the EP02 problem, similar performance for resolving the density and pressure profiles is obtained for the various TENO schemes and the CWENO5 scheme, as shown in Fig.~\ref{Fig:2d_explosion_ep02_01},~\ref{Fig:2d_explosion_ep02_02}, and~\ref{Fig:2d_explosion_ep02_03}. However, the WENO5 scheme generates spurious numerical oscillations near the discontinuities, which can be observed from Fig.~\ref{Fig:2d_explosion_ep02_02}(a), Fig.~\ref{Fig:2d_explosion_ep02_03}(c), and Fig.~\ref{Fig:2d_explosion_ep02_03}(d).

\subsection{2D Kelvin-Helmholtz instability (KHI)}
    The Kelvin-Helmholtz instability with single mode perturbation is considered \cite{san2015evaluation}\cite{ryu2000magnetohydrodynamic}. The computational domain is $[-0.5, 0.5] \times [-0.5, 0.5]$ and the initial conditions are given as
    \begin{equation}
        (\rho ,{u}) = \left\{ {\begin{array}{*{20}{c}}
        {(2, - 0.5),}&{\left| y \right| \le 0.25,}\\
        {(1,0.5),}&{\text{otherwise,}} \\
        \end{array}} \right. \\
    \end{equation}
    and
    \begin{equation}
        ({v},p) = (0.01\sin (2\pi x),2.5).
    \end{equation}
    $\gamma = 1.4$ and periodic boundary conditions are enforced at the boundaries of the computational domain. The final simulation time is $t=1$.

\begin{figure}[H]
        \centering
        \includegraphics[width=1.0\textwidth,trim=2 2 2 2,clip]{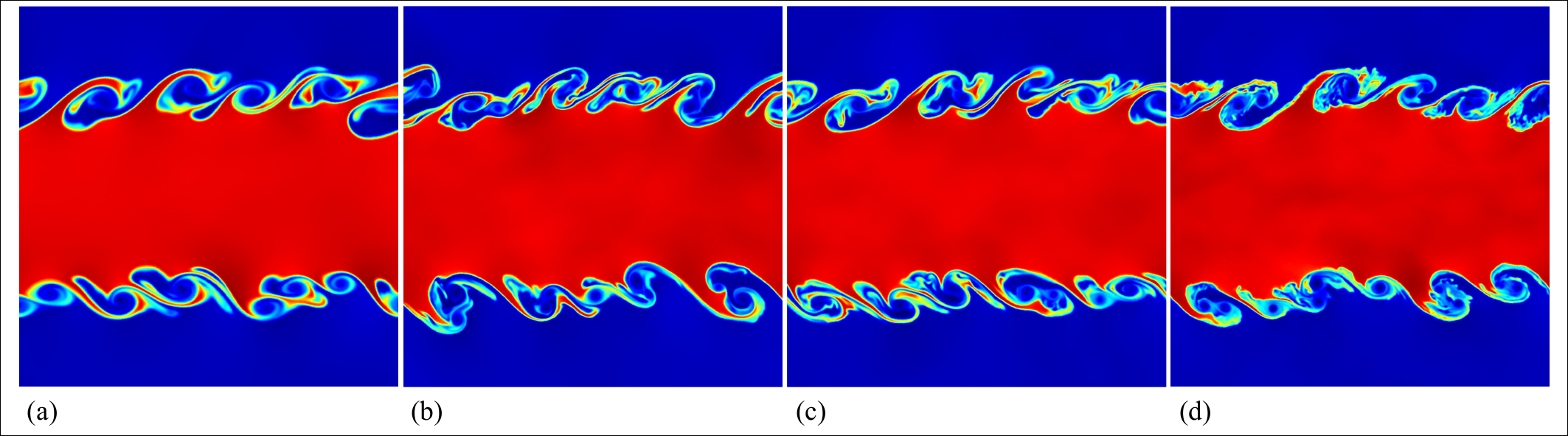}
        \caption{2D Kelvin-Helmholtz instability: density distributions at $t=1$ from the (a) TENO3, (b) TENO4, (c) TENO5 and (d) TENO6 schemes. The uniformly distributed triangular mesh with an effective resolution $h\approx1/500$ is adopted.}
        \label{Fig:2d_khi_01}
    \end{figure}

\begin{figure}[H]
        \centering
        \includegraphics[width=1.0\textwidth,trim=2 2 2 2,clip]{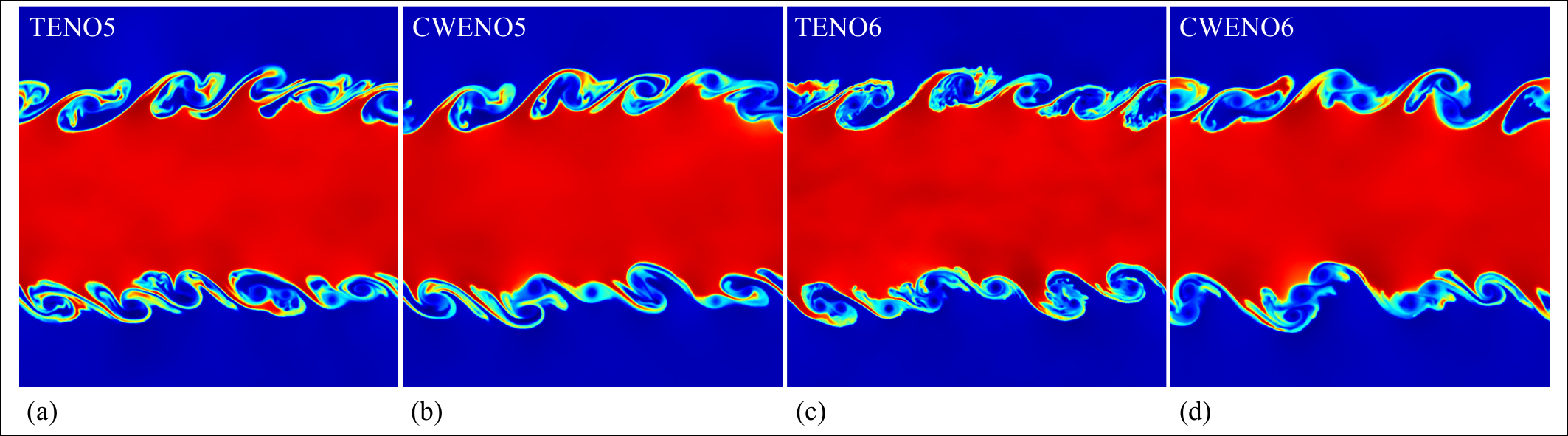}
        \caption{2D Kelvin-Helmholtz instability: density distributions at $t=1$ from the (a) TENO5, (b) CWENO5, (c) TENO6 and (d) CWENO6 schemes. The uniformly distributed triangular mesh with an effective resolution $h\approx1/500$ is adopted.}
        \label{Fig:2d_khi_02}
    \end{figure}

\begin{figure}[H]
        \centering
        \includegraphics[width=1.0\textwidth]{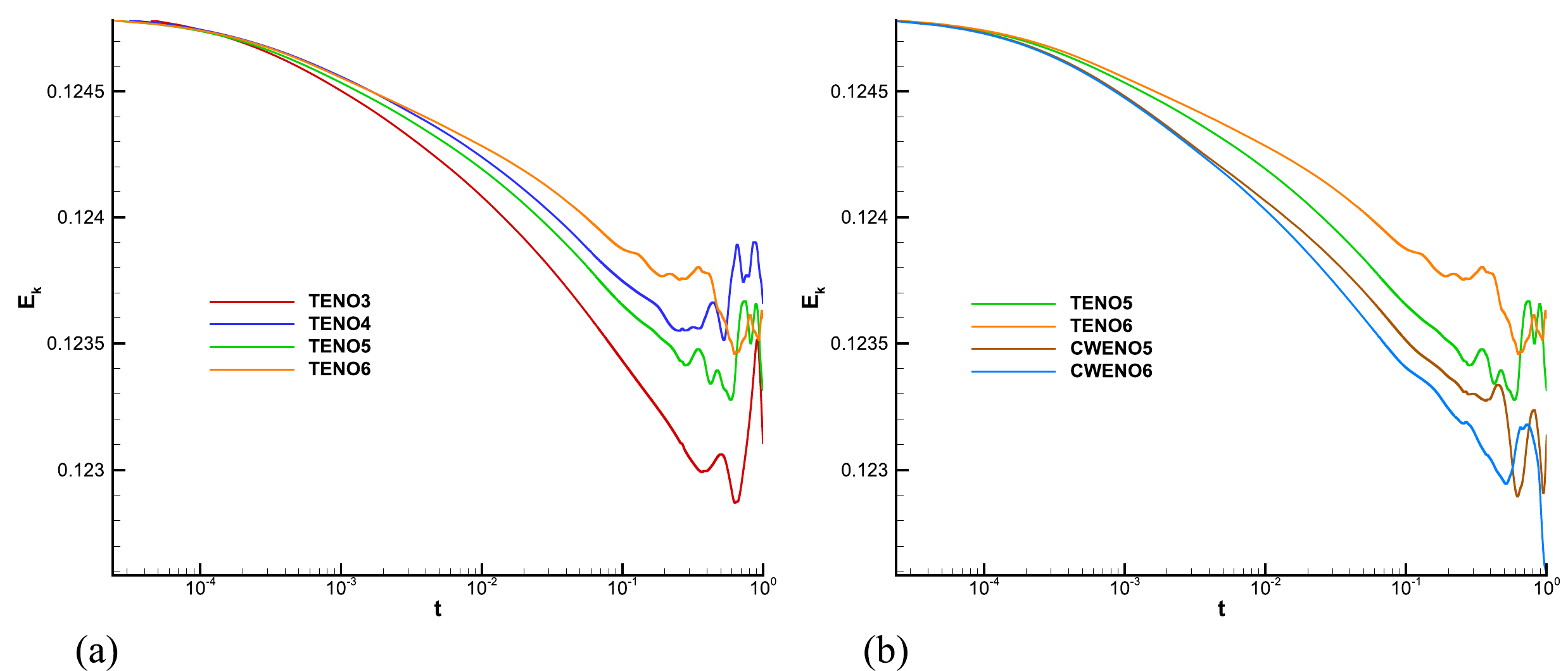}
        \caption{2D Kelvin-Helmholtz instability: density profiles from (a) TENO schemes of different orders and comparisons between (b) TENO5/6 and CWENO5/6. The uniformly distributed triangular mesh with an effective resolution $h\approx1/500$ is adopted.}
        \label{Fig:2d_khi_03}
    \end{figure}

    The density distributions computed by various TENO, WENO and CWENO schemes are given by Fig.~\ref{Fig:2d_khi_01},~\ref{Fig:2d_khi_02}, and~\ref{Fig:2d_khi_03}. The third-order TENO3 scheme is much more dissipative than other TENO schemes of higher order and smears the small-scale vortical structures significantly. As shown in Fig.~\ref{Fig:2d_khi_03}, the TENO5 and TENO6 schemes perform much better than the CWENO5 and CWENO6 schemes in terms of preserving the kinetic energy.

\subsection{2D double Mach reflection of a strong shock (DMR)}
\begin{figure}[H]
        \centering
        \includegraphics[width=1.0\textwidth]{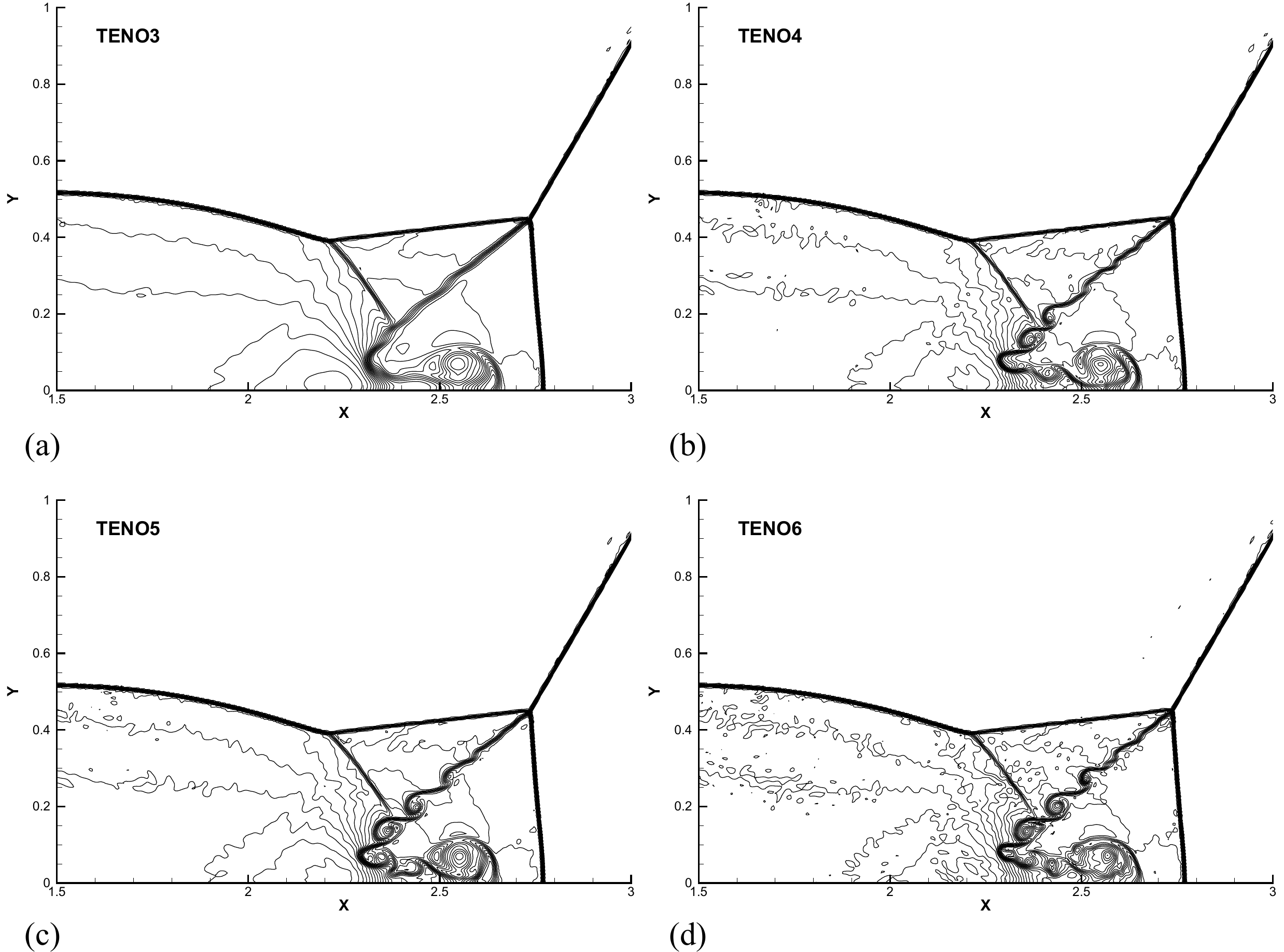}
        \caption{2D double Mach reflection of a strong shock: results from the (a) TENO3, (b) TENO4, (c) TENO5 and (d) TENO6 schemes. The figures are drawn with 43 density contours between 1.887 and 20.9. The mesh has 369,794 uniform triangular elements and $h\approx1/200$.}
        \label{Fig:dmr_02_mesh4_tenos}
    \end{figure}

\begin{figure}[H]
        \centering
        \includegraphics[width=1.0\textwidth]{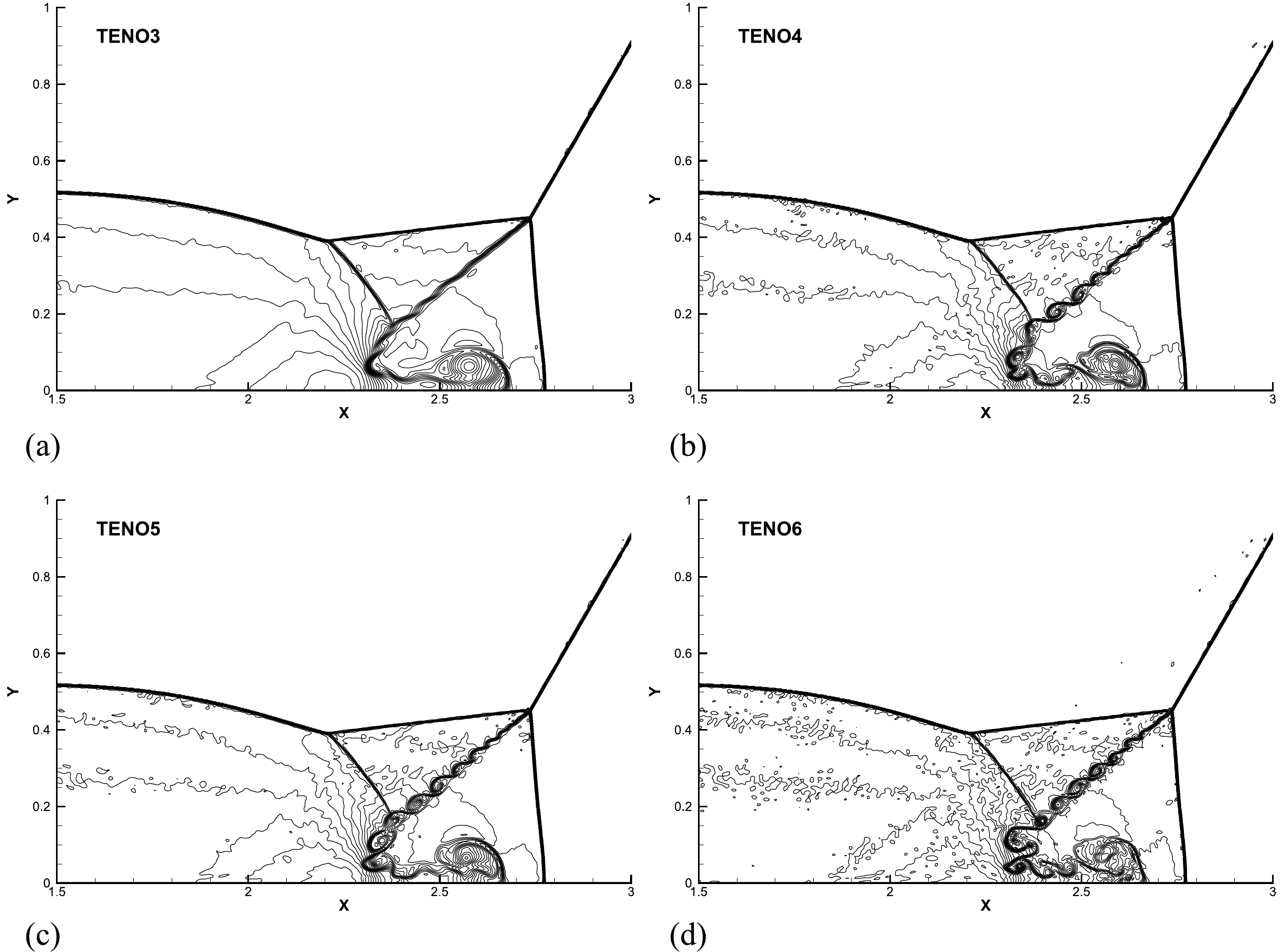}
        \caption{2D double Mach reflection of a strong shock: results from the (a) TENO3, (b) TENO4, (c) TENO5 and (d) TENO6 schemes. The figures are drawn with 43 density contours between 1.887 and 20.9. The mesh has 752,708 uniform triangular elements and $h\approx1/280$.}
        \label{Fig:dmr_06_mesh5_tenos}
    \end{figure}

\begin{figure}[H]
        \centering
        \includegraphics[width=1.0\textwidth,trim=2 2 2 2,clip]{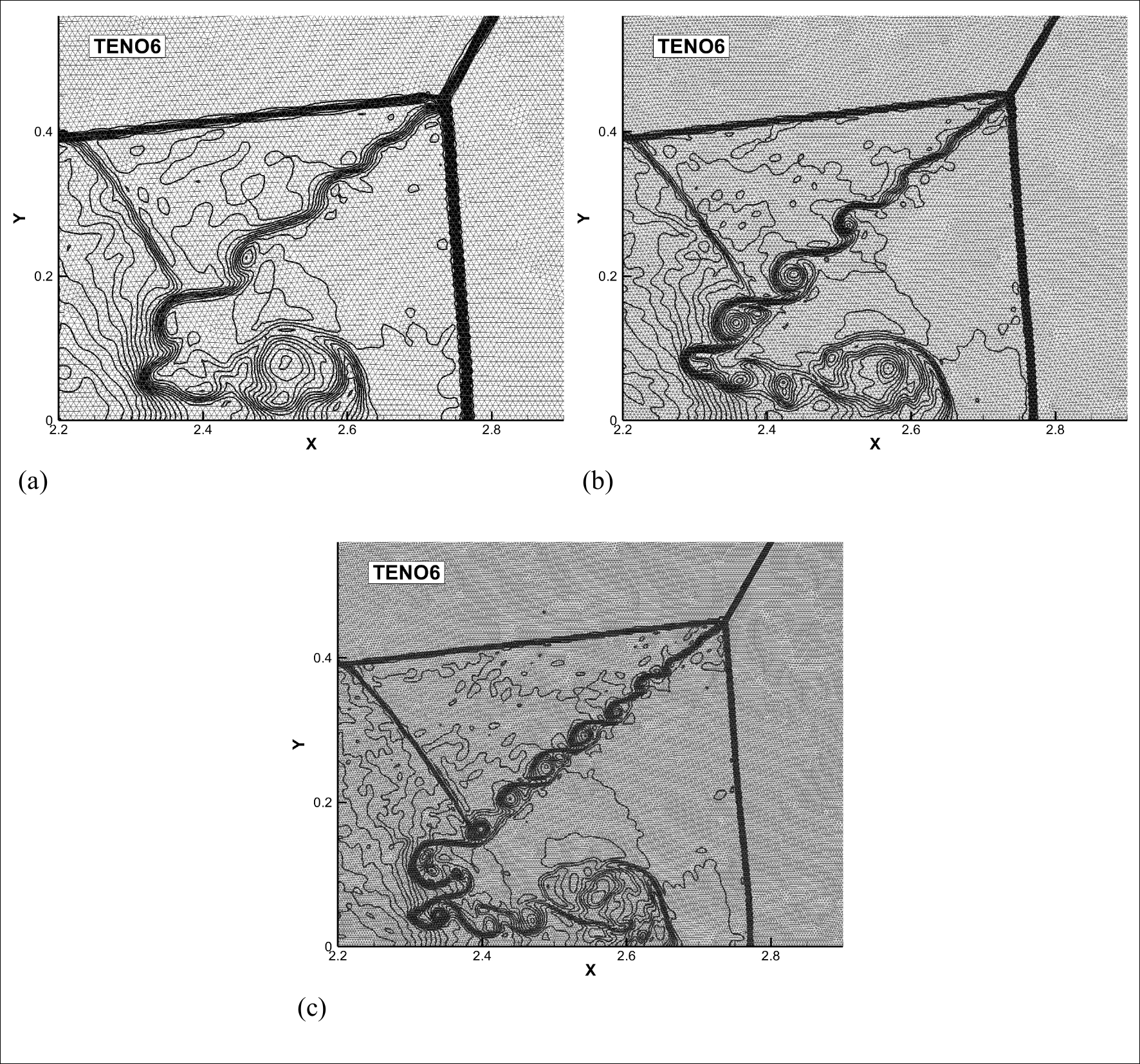}
        \caption{2D double Mach reflection of a strong shock: zoomed-in results from the TENO6 scheme with (a) 151,850, (b) 369,794 and (c) 752,708 uniform triangular elements. The figures are drawn with 43 density contours between 1.887 and 20.9. Also plotted are the mesh topologies.}
        \label{Fig:dmr_05_teno6_with_mesh}
    \end{figure}

\begin{figure}[H]
        \centering
        \includegraphics[width=1.0\textwidth]{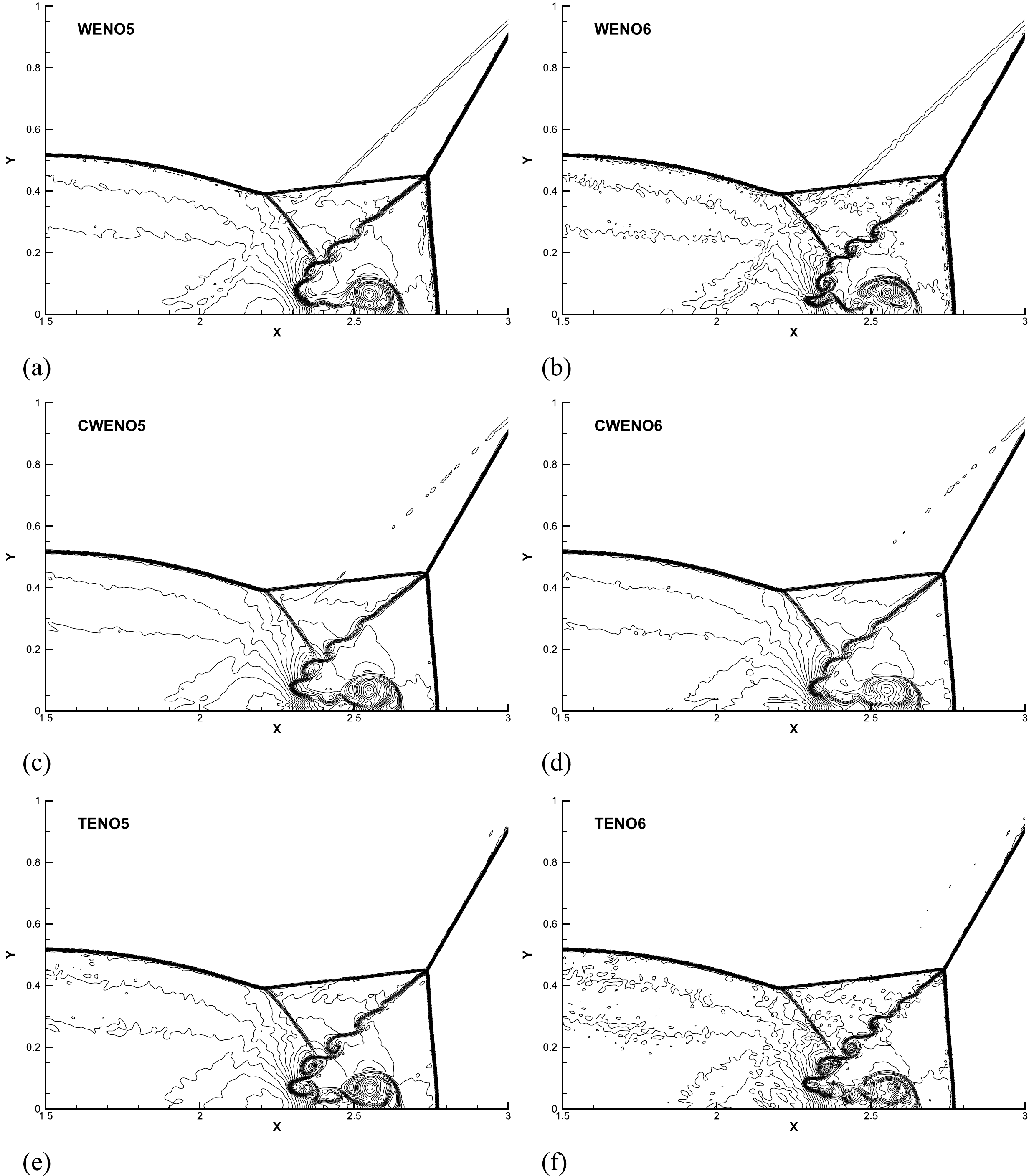}
        \caption{2D double Mach reflection of a strong shock: results from the (a) WENO5, (b) WENO6, (c) CWENO5, (d) CWENO6, (e) TENO5, (f) TENO6 schemes. The figures are drawn with 43 density contours between 1.887 and 20.9. The mesh has 369,794 uniform triangular elements and $h\approx1/200$. Note that both the WENO5 and WENO6 schemes fail at this resolution without the additional limiter in \cite{tsoutsanis2018extended}.}
        \label{Fig:dmr_04_mesh4_teno_weno_cweno}
    \end{figure}

\begin{figure}[H]
        \centering
        \includegraphics[width=1.0\textwidth]{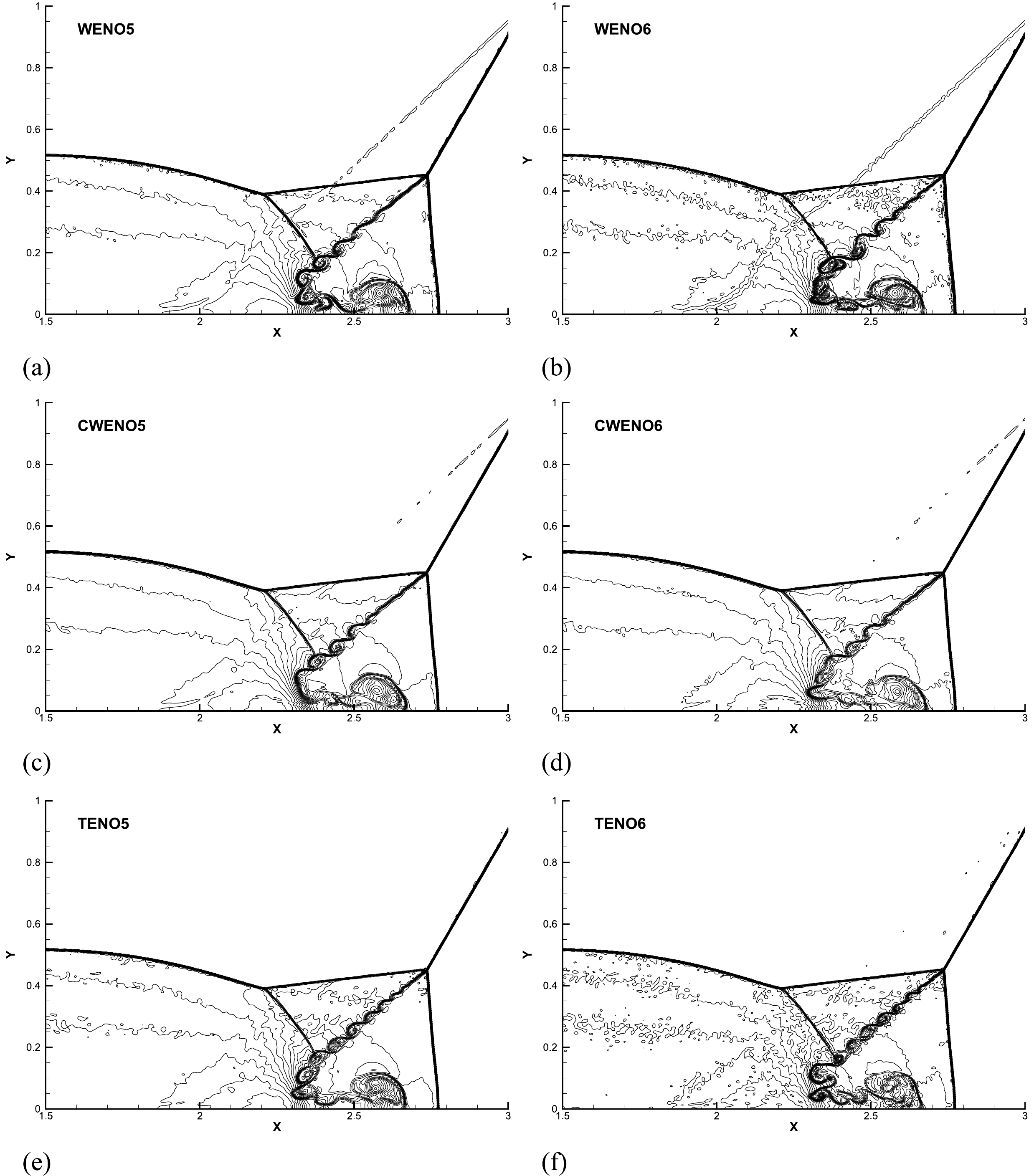}
        \caption{2D double Mach reflection of a strong shock: results from the (a) WENO5, (b) WENO6, (c) CWENO5, (d) CWENO6, (e) TENO5, (f) TENO6 schemes. The figures are drawn with 43 density contours between 1.887 and 20.9. The mesh has 752,708 uniform triangular elements and $h\approx1/280$. Note that the WENO6 scheme fails at this resolution without the additional limiter in \cite{tsoutsanis2018extended}.}
        \label{Fig:dmr_07_mesh5_teno_weno_cweno}
    \end{figure}

{\color{black}The initial condition is \cite{Woodward1984}
    \begin{equation}
        \label{eq:DMR}
        (\rho ,u,v,p) = \left\{ {\begin{array}{*{20}{c}}
        {(1.4,0,0,1) ,}&{\text{if } y < 1.732(x - 0.1667) ,}\\
        {(8,7.145, - 4.125,116.8333) , }&{\text{otherwise }.}
        \end{array}} \right.
    \end{equation}
    The computational domain is $[0,4] \times [0,1]$ and the simulation end time is $t = 0.2$. Initially, a right-moving Mach 10 shock wave is placed at $x = 0.1667$ with an incident angle of $60^{\circ}$ to the x-axis. The post-shock condition is imposed from $x = 0$ to $x = 0.1667$ whereas a reflecting wall condition is enforced from $x = 0.1667$ to $x = 4$ at the bottom. For the top boundary condition, the fluid variables are defined to exactly describe the evolution of the Mach 10 shock wave. The inflow and outflow condition are imposed for the left and right sides of the computational domain.
    
    The numerical results from TENO schemes of various orders and the corresponding CWENO and WENO schemes are given in Fig.~\ref{Fig:dmr_02_mesh4_tenos},~\ref{Fig:dmr_06_mesh5_tenos},~\ref{Fig:dmr_05_teno6_with_mesh},~\ref{Fig:dmr_04_mesh4_teno_weno_cweno}, and~\ref{Fig:dmr_07_mesh5_teno_weno_cweno}. As shown in~\ref{Fig:dmr_02_mesh4_tenos} and~\ref{Fig:dmr_06_mesh5_tenos}, with the same mesh resolution, the TENO schemes with higher-order resolve the small-scale features much better with less numerical dissipation. The zoomed-in view of the density distributions in Fig.~\ref{Fig:dmr_05_teno6_with_mesh} shows that no spurious numerical oscillations present in the vicinity of the discontinuities. As given by Fig.~\ref{Fig:dmr_04_mesh4_teno_weno_cweno} and~\ref{Fig:dmr_07_mesh5_teno_weno_cweno}, when compared with the WENO and CWENO schemes of same orders, the TENO schemes feature the least numerical dissipation, followed by the WENO schemes and subsequently by the CWENO schemes. It is worth noting that, while both the TENO and CWENO schemes are robust for all the considered mesh resolutions, the WENO5 and WENO6 schemes typically fail the simulations without an additional limiter \cite{tsoutsanis2018extended}, which helps preserve the positivity of density and pressure.
    
    Moreover, when compared to the newly proposed fifth-order WENO5-AO scheme, the present TENO5 scheme performs better in capturing the roll-up of the Mach stem with half resolution in each coordinate direction, see their Fig.~4 of \cite{balsara2020efficient}. The result from the present TENO5 scheme at the resolution of $h\approx1/200$ is comparable to that from the multi-resolution WENO5 scheme at the the resolution of $h\approx1/320$ in \cite{zhu2019new}.}

\subsection{{\color{black}3D explosion problem}}
{\color{black}We consider a 3D explosion problem to validate the robustness of high-order TENO schemes on unstructured meshes. The initial condition is given as
    \begin{equation}
        \label{eq:bw_ep1}
        (\rho,u,v,w,p)(\textbf{r},0)=
        \begin{cases}
        (1.0,0,0,0,1.0), & \text{if $ \Vert \textbf{r} \Vert \leqslant 0.5$},\\
        (0.125,0,0,0,0.1), & \text{otherwise}.
        \end{cases}
    \end{equation}
The final simulation time is $t = 0.25$. A uniform tetrahedron mesh with edge length $h\approx1/30$ is adopted.}

\begin{figure}[H]
        \centering
        \includegraphics[width=0.9\textwidth]{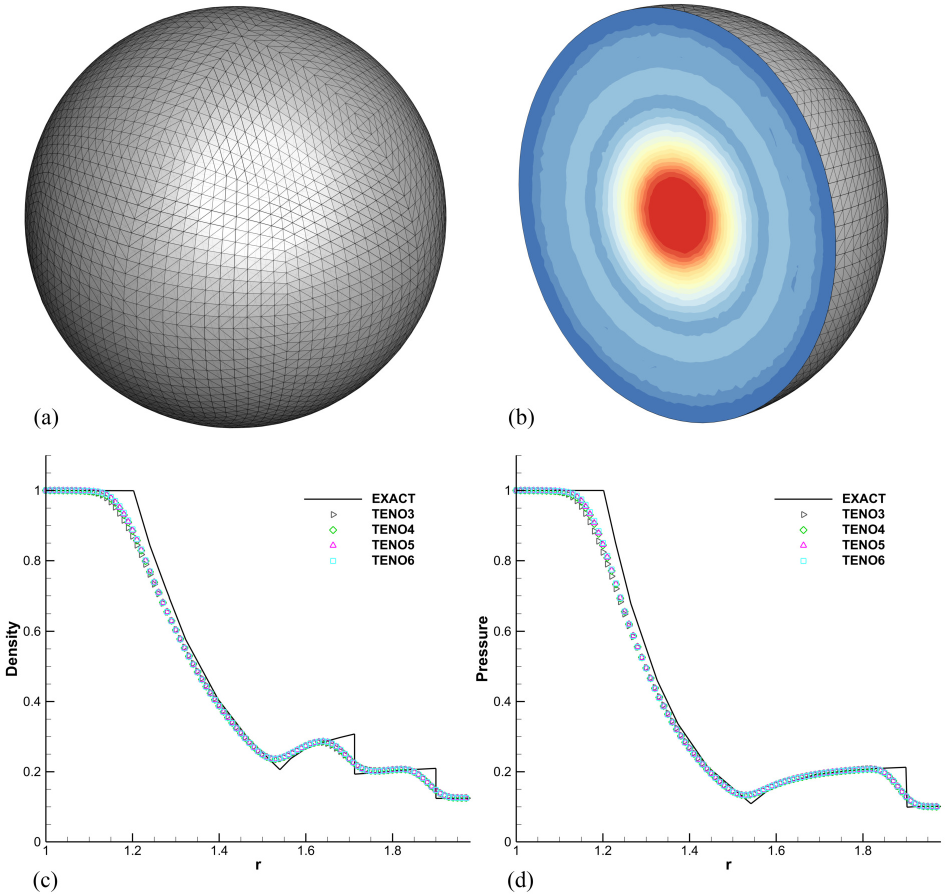}
        \caption{{\color{black}3D explosion problem: (a) and (b) denote the deployed uniform tetrahedron mesh with edge length $h\approx1/30$ and the density contour from the TENO4 scheme; (c) and (d) denote the comparisons of the density and pressure profiles, respectively.}}
        \label{Fig:explosion}
    \end{figure}

{\color{black}As shown in Fig.~\ref{Fig:explosion}, no obvious numerical artifacts can be identified from the comparisons of the density and pressure profiles. Quantitatively, the resolved solutions from all the considered TENO schemes agree with the analytical solutions well.}

    \color{black}

\subsection{Efficiency analysis}

    Table \ref{Tab:runtime_analysis} lists the statistics regarding the computing time of different schemes, where $N_{iter}$ denotes the number of iterations to complete the calculation, $t_{total}$ the total computing time in seconds and $t_{recon}$ the computing time for the high-order reconstruction. Upon further normalization by $N_{iter}$ and $N_{e}$, we obtain $t_{norm}$ (second/iteration/element) from $t_{recon}$. Two cases are considered for the analysis, i.e. the 2D DMR and the 2D KHI. The total number of mesh elements is 369,794 and 577,416, respectively. For DMR, we measure the total computing time at $t=2$ and for KHI, the total computing time is measured at $t=0.1$. 

    By increasing the accuracy order of TENO schemes from third- to sixth-order, the normalized cost increases by roughly $65\%$ for the same simulation. Overall speaking, the TENO and CWENO schemes with the same accuracy order feature a similar computational efficiency for both cases.
    \begin{table}[H]
    \centering
    \begin{threeparttable}
        \caption{Computing-time statistics for TENO and CWENO schemes}
        \scriptsize
        \label{Tab:runtime_analysis}
        \newcommand{\tabincell}[2]{\begin{tabular}{@{}#1@{}}#2\end{tabular}}
        \begin{tabular}{>{\centering\arraybackslash}m{2.0cm}
                        >{\centering\arraybackslash}m{2.0cm}
                        >{\centering\arraybackslash}m{2.0cm}
                        >{\centering\arraybackslash}m{2.0cm}
                        >{\centering\arraybackslash}m{2.0cm}
                        >{\centering\arraybackslash}m{2.0cm}}
        \hline
             & scheme & $N_{iter}$ & $t_{total}$ & $t_{recon}$ & $t_{norm}$ \\ \hline
         DMR & TENO3  & 16007  & 9468.00  & 8079.77  & 1.36E-06 \\
             & TENO4  & 16110  & 12194.47 & 9895.88  & 1.66E-06 \\
             & TENO5  & 16081  & 13990.78 & 11674.00 & 1.96E-06 \\
             & TENO6  & 16152  & 16262.73 & 13358.49 & 2.24E-06 \\
             & CWENO5 & 16045  & 14126.61 & 11738.53 & 1.98E-06 \\
             & CWENO6 & 16019  & 15472.95 & 13016.93 & 2.20E-06 \\
        \hline
         KHI & TENO3  & 3580   & 3298.40  & 2733.38  & 1.32E-06 \\
             & TENO4  & 3580   & 4201.04  & 3425.22  & 1.66E-06 \\
             & TENO5  & 3579   & 4587.06  & 3935.91  & 1.90E-06 \\
             & TENO6  & 3579   & 5399.16  & 4540.90  & 2.20E-06 \\
             & CWENO5 & 3578   & 4626.61  & 3887.51  & 1.88E-06 \\
             & CWENO6 & 3579   & 5414.79  & 4491.87  & 2.17E-06 \\
        \hline
        \end{tabular}
    \begin{tablenotes}
      \small
      \item $\ast$ \textcolor{black}{Regarding the hardware configurations, the CPU used for the tests is Intel Xeon E5-2690. The Intel Fortran Compiler version 19 and the Intel MPI Library 2019 are employed for the compiling with the compiler flags \textit{-i4 -r8 -ipo -march=core-avx2 -mtune=core-avx2 -O3 -fp-model precise -zero -qopenmp  -qopenmp-link=static}. For both the DMR and KHI simulations, 56 parallel threads are employed, i.e., 8 MPI tasks with 7 OpenMP threads per MPI rank.}
    \end{tablenotes}
    \end{threeparttable}
    \end{table}

\section{Conclusions}

In this paper, the high-order TENO scheme, originally proposed for structured meshes, is for the first time extended to unstructured meshes. A set of TENO schemes from third- to sixth-order accuracy is developed and validated. The concluding remarks are given as follows.

    \begin{itemize}
        \item The new candidate stencils include one large central-biased stencil and several small directional stencils. The targeted high-order reconstruction scheme is constructed on the large stencil while the third-order schemes are constructed on the small stencils. To achieve the same accuracy order in smooth regions, the present stencil arrangement leads to narrower full stencil width when compared to those deployed in  \cite{tsoutsanis2011weno}\cite{titarev2010weno}\cite{dumbser2007arbitrary}\cite{dumbser2007quadrature}.
        
        \item Following a strong scale separation procedure, a novel ENO-like stencil selection strategy is tailored for unstructured meshes. Such a TENO weighting strategy ensures that the high-order accuracy is restored by adopting the candidate reconstruction from the large stencil and the sharp shock-capturing capability is retained by selecting the candidate reconstruction from the smooth small stencils. As analyzed in previous work \cite{fu2016family}\cite{fu2018new}\cite{fu2019low}, one dissipation source of WENO family schemes comes from the smooth weighting of the contributions from candidate stencils, which leads to the early triggering of the nonlinear adaptations in the low-wavenumber regimes. Similar to those TENO schemes for Cartesian meshes, the low-dissipation property is inherited in the proposed schemes due to the sharp stencil selection, which either applies the candidate stencil in an optimal way or abandons it completely when crossed by genuine discontinuities. 
        
        \item The new framework allows for arbitrarily high-order TENO reconstructions. For conceptual investigations, the TENO schemes from third- to sixth-order accuracy are constructed and the built-in parameters are explicitly given.
        
        \item A set of challenging benchmark simulations has been conducted with the default set of parameters. Numerical results reveal that the proposed TENO schemes are robust for highly compressible gas dynamic simulations with low numerical dissipation and sharp discontinuity-capturing capability.
        
        \item Compared to the WENO schemes with equal-size candidate stencils, the overall stencil width of the present TENO schemes is much more compact. Compared to the CWENO schemes with the same accuracy order, the proposed TENO schemes are much less dissipative while featuring the sharp shock-capturing capability. Note that, the present investigation also reveals that the WENO schemes with equal-size candidate stencils are less robust when deployed to strong discontinuities (as shown in the DMR simulations), which denotes another flaw in addition to the complex reconstruction procedure.
    \end{itemize}

    Considering the good performance of the proposed TENO framework and the flexibility for further extension, the future work will focus on the deployment of the present schemes to more complicated flows, e.g. the chemical reacting flows, the MHD flows, and the external aerodynamics with realistic geometries.
    
    {\color{black}Another potential research direction is to develop TENO schemes with lower-order more compact directional stencils, which may provide additional benefits, e.g., further reducing the communication overheads for large-scale parallel computations and preserving better numerical robustness in regions with strong gradient or poor grid quality. The performance as well as the convergence studies for this type of TENO schemes will be investigated and reported separately in a forthcoming paper.}

\section*{Declaration of Competing Interest}

    The authors declare that they have no known competing financial interests or personal relationships that could have appeared to influence the work reported in this paper.
\section*{Data availability}
The data that support the findings of this study are available on request from the corresponding author, Lin Fu.

\section*{Acknowledgements}

L.F. acknowledges the fund from Guangdong Basic and Applied Basic Research Foundation (No. 2022A1515011779), the fund from Shenzhen Municipal Central Government Guides Local Science and Technology Development Special Funds Funded Projects (No. 2021Szvup138) and the fund from CORE as a joint research center for ocean research between QNLM and HKUST.

\bibliographystyle{elsarticle-num}


\end{document}